\documentclass{article}

\usepackage[margin=1in]{geometry}
\usepackage{amsmath,amssymb,amsfonts,graphicx,amsthm,nicefrac,mathtools,bm}
\usepackage{hyperref}
\usepackage[numbers]{natbib}
\usepackage[english]{babel}
\usepackage{times}
\usepackage[T1]{fontenc}
\usepackage{bbm}
\usepackage{color}
\usepackage{xspace}
\usepackage{rotating,multirow}
\usepackage{caption}
\usepackage{subcaption}

\newtheorem{theorem}{Theorem}

\newtheorem{lemma}{Lemma}

\newtheorem{corollary}{Corollary}

\newtheorem{assumption}{Assumption}
\theoremstyle{definition}
\newtheorem{definition}{Definition}
\theoremstyle{remark}

\makeatletter
\newtheorem*{rep@theorem}{\rep@title}
\newcommand{\newreptheorem}[2]{%
\newenvironment{rep#1}[1]{%
 \def\rep@title{#2 \ref{##1}}%
 \begin{rep@theorem}}%
 {\end{rep@theorem}}}
\makeatother
\newreptheorem{theorem}{Theorem}
\newreptheorem{lemma}{Lemma}
\newreptheorem{corollary}{Corollary}
\newreptheorem{proposition}{Proposition}

\newcommand{\assumpref}[1]{Assumption~\ref{assump:#1}}
\newcommand{\assumpsref}[1]{Assumptions~\ref{assump:#1}}
\newcommand{\assumpssref}[1]{\ref{assump:#1}}
\newcommand{\figref}[1]{Figure~\ref{fig:#1}}

\newcommand{\secref}[1]{Section~\ref{sec:#1}}
\newcommand{\appref}[1]{Appendix~\ref{app:#1}}

\newcommand{\defref}[1]{Definition~\ref{def:#1}}

\newcommand{\lemref}[1]{Lemma~\ref{lem:#1}}
\newcommand{\lemsref}[1]{Lemmas~\ref{lem:#1}}
\newcommand{\lemssref}[1]{\ref{lem:#1}}

\newcommand{\thmref}[1]{Theorem~\ref{thm:#1}}

\newcommand{\eqnref}[1]{\eqref{eqn:#1}}

\DeclareMathOperator{\rank}{rank}

\DeclareMathOperator*{\argmin}{arg\,min}

\newcommand{\inner}[2]{\langle{#1},{#2}\rangle} 
\newcommand{\norm}[1]{\lVert{#1}\rVert}
\newcommand{\Norm}[1]{\left\lVert{#1}\right\rVert}

\renewcommand{\O}[1]{\mathcal{O}\left({#1}\right)}
\def\R{\mathbb{R}}

\newcommand{\ident}{\mathbf{I}}

\newcommand{\iidsim}{\stackrel{\mathrm{iid}}{\sim}}

\newcommand{\ignore}[1]{}
\renewcommand{\Pr}[2]{\mathcal{P}_{{#1}}\left({#2}\right)} 

\newcommand{\pr}[1]{\mathcal{P}_{{#1}}}

\newcommand{\eps}{\epsilon}

\usepackage{fancyhdr}
\newcommand{\thedate}{\today}
\newcommand{\theauthor}{Wooseok Ha and Rina Foygel Barber}
\newcommand{\thetitle}{Alternating minimization and alternating  descent\\
 over nonconvex sets}

\date{\thedate}
\author{\theauthor}
\title{\thetitle}

\fancypagestyle{plain}{%
  \fancyhf{}%
  \lhead{\thepage}
}

\usepackage[mmddyy]{datetime}

\newcommand{\loss}{\mathcal{L}}

\newcommand{\errx}{\epsilon_{x}}
\newcommand{\erry}{\epsilon_{y}}

\newcommand{\biginner}[2]{\left\langle{#1},{#2}\right\rangle} 

\newcommand{\opnorm}[1]{\norm{{#1}}_{\textnormal{op}}}
\newcommand{\nucnorm}[1]{\norm{{#1}}_{\textnormal{nuc}}}
\newcommand{\opNorm}[1]{\Norm{{#1}}_{\textnormal{op}}}
\newcommand{\fronorm}[1]{\norm{{#1}}_{\textnormal{F}}}

\newcommand{\xh}{\widehat{x}}
\newcommand{\yh}{\widehat{y}}

\newcommand{\xt}{\widetilde{x}}

\newcommand{\Uset}{\mathcal{U}}

\newcommand{\Xset}{\mathcal{X}}
\newcommand{\Yset}{\mathcal{Y}}

\newcommand{\Zset}{\mathcal{Z}}

\newcommand{\vect}[1]{\textnormal{vec}\left({#1}\right)}
\newcommand{\ufun}{\mathsf{g}}

\newcommand{\Xs}{X^\star}
\newcommand{\Ys}{Y^\star}
\newcommand{\Xh}{\widehat{X}}
\newcommand{\Yh}{\widehat{Y}}
\newcommand{\Us}{U^\star}

\newcommand{\Aoper}{\mathcal{A}}

\newcommand{\spike}{\alpha_{\mathsf{sp}}}
\newcommand{\Sigmas}{\Sigma^\star}
\newcommand{\Sigmah}{\widehat{\Sigma}}

\newcommand{\Thetas}{\Theta^\star}
\newcommand{\Thetah}{\widehat{\Theta}}
\newcommand{\Thetaset}{\mathcal{Q}}
\newcommand{\lammax}{\lambda_{\textnormal{max}}}
\newcommand{\lammin}{\lambda_{\textnormal{min}}}

\newcommand{\normx}[1]{\norm{{#1}}_{x}}
\newcommand{\normy}[1]{\norm{{#1}}_{y}}

\newcommand{\ball}{\mathbb{B}}

\newcommand{\diver}[2]{{D}^2(#1;#2)}
\newcommand{\sdiver}[2]{{D}(#1;#2)}

\begin{document}

\maketitle

\begin{abstract}
We analyze the performance of alternating minimization for loss functions optimized over two variables, where each variable may be restricted to lie in some potentially
nonconvex constraint set.
This type of setting arises naturally in high-dimensional statistics and signal processing, where the variables often reflect different structures or components within the signals being considered. Our analysis relies on the notion of local concavity coefficients, which has been proposed in \citet{barber2017gradient} to measure and quantify the concavity of a general nonconvex set. Our results further reveal important distinctions between alternating and non-alternating methods. Since computing the alternating minimization steps may not be tractable for some problems, we also consider an inexact version of the algorithm and provide a set of sufficient conditions to ensure fast convergence of the inexact algorithms. We demonstrate our framework on several
examples, including low rank + sparse decomposition and multitask regression, and provide numerical experiments to validate our theoretical results.
\end{abstract}

\section{Introduction}\label{sec:intro}

Many methods in modern statistics use structured constraints to improve signal recovery in a high-dimensional setting. 
Common constraints in the high-dimensional statistics literature include sparsity, requiring that a signal consists of mostly zero values; variants on sparsity, such as total variation sparsity, requiring that the first-order differences of a signal are locally constant; and low rank, where the signal is a matrix expressed as the sum of only a few linear factors. In many settings, multiple structures may be present simultaneously in the data,
in which case we may need to optimize a function over several variables,  which are each believed to exhibit some latent
structure---for instance, a low-rank term and a sparse term.

Much of the literature in this area focuses on convex relaxations of these structured constraints, such as the $\ell_1$ norm (as a convex approximation to sparsity). Working
with a convex penalty or convex constraint, as a proxy for the nonconvex structure of the variable(s) of interest, allows for easier optimization from both a theoretical and a practical
point of view. In recent years, however, attention has turned to nonconvex optimization problems, aiming to avoid the loss of accuracy that is often the cost of taking a convex
relaxation. 

In this work, we consider the problem of optimizing over two variables, one or both of which is constrained to lie in some potentially nonconvex set:
\begin{equation*}\label{eqn:altmin_program}(\xh,\yh) = \argmin\{\loss(x,y):x\in\Xset,y\in\Yset\},\end{equation*}
where $\Xset\subset\R^{d_x}$ and $\Yset\subset\R^{d_y}$ reflect our beliefs or desired properties for the $x$ and $y$ variables,
while $\loss$ is the target function to minimize (for example, a negative log-likelihood, in which case we are searching for the constrained maximum likelihood
estimator).

Our aim is to study the convergence behavior of two methods for this problem. First, we consider the {\em alternating minimization}
method, where we iterate the steps
\[\begin{cases} 
\text{Fix $y$, and choose $x \in \Xset$ to minimize the function $x\mapsto\loss(x,y)$};\\
\text{Fix $x$, and choose $y \in \Yset$ to minimize the function $y\mapsto\loss(x,y)$}.
\end{cases}\]
This type of method can be practical in scenarios where the loss function is relatively simple to minimize when viewed as a function of either $x$ or $y$ only---for instance,
in multivariate regression, where $x$ represents the coefficients and $y$ represents the covariance structure.
In other settings, even the marginal minimization steps are expensive to calculate, but we can instead consider approximating each one with gradient descent.
The resulting alternating gradient descent algorithm iterates the following steps:
\[\begin{cases} 
\text{Fix $y$, and update $x$ by taking (one or more) gradient descent steps on the function $x\mapsto\loss(x,y)$};\\
\text{Fix $x$, and update $y$ by taking (one or more) gradient descent steps on the function $y\mapsto\loss(x,y)$}.
\end{cases}\]

Our main results derive conditions under which both of these algorithms converge linearly, 
under certain assumptions:
\begin{itemize}
\item Loss function:
The assumptions on the loss function $\loss(x,y)$ are familiar in the high-dimensional statistics literature, namely,
restricted strong convexity (RSC) and restricted smoothness (RSM) assumptions, which essentially require that $\loss(x,y)$
behaves like a smooth and strongly convex function when restricted to the constrained domains $\Xset$ and $\Yset$.
\item Nonconvex constraints: 
For the constraint sets $\Xset$ and $\Yset$, we work in the framework established in \citet{barber2017gradient},
requiring bounded {\em local concavity coefficients} for each set (see \secref{nonconvex_constraint} for further details). This geometric condition is a natural relaxation of convexity
and is  satisfied by many commonly used
nonconvex constraints, such as a low-rank constraint.
\item Initialization: In order to ensure convergence to a global minimum, we need to assume that the initialization point
is sufficiently close so that the convexity of the loss function $\loss(x,y)$ is sufficient to outweigh local concavity in the constraint
sets. Details are given in our results below.
\end{itemize}
To demonstrate the utility of our results, we also consider a range of specific examples with rank-constrained variables,
including factor models, multivariate regression, and robust PCA.

\subsection{Related work}\label{sec:relatedwork}

Alternating minimization is a classical topic in the optimization literature, and a large body of research has been devoted to understanding the method under various settings. On the other hand, nonconvex constraints have recently received a lot of attention from community, including many results treating sparsity-constrained or rank-constrained problems specifically. Our work combines both settings, and thus is naturally related to many existing works. Here, we summarize some of the key recent results, and describe how they relate to our contributions; for brevity, we only focus on the papers most relevant to our work.

\paragraph{Nonconvex constraints on a single variable} 

The past few years have witnessed extensive results on the optimization problem over a nonconvex set. For instance, \citet{jain2014iterative} consider the problem of minimizing a loss function over a sparsity constraint or a rank constraint, and show that the iterative hard thresholding method can achieve global convergence to a target point, as long as the sparsity or rank of the target is far smaller than the given threshold of the constraint. For a rank constraint specifically, \citet{grussler2016low} study the duality-gap for a certain class of problems and present some situations under which there exists no duality-gap; therefore, solutions to the original problem and its convex relaxation (largest convex minorizer) must coincide under these situations. Turning to a more general setting, \citet{oymak2015sharp} study projected gradient descent scheme for least squares objective when constraining to a nonconvex regularizer set, establishing linear convergence of the algorithm from any initialization point. In particular, by introducing a descent cone at the target point, the authors  characterize the convergence rate in terms of the singular values of the Hessian matrix restricted to this cone.  The work of \citet{barber2017gradient} takes a different approach and instead develops a way of measuring local concavity of the constraint set at any given point. This measure of local concavity
is then used to  analyze the local convergence of projected gradient descent. 

In the setting of two variables, $x$ and $y$, as considered in this paper, \citet{barber2017gradient}'s approach can be applied by performing gradient
descent on the joint variable $(x,y)$, constraining to the space $\Xset\times \Yset$---that is, any problem with multiple variables can of course  be reformulated
as a single variable problem. In this work, however, we find that separating the variables and alternating their updates can provide substantial benefits, both
theoretically and empirically. We discuss these issues in detail in \secref{altmin}.


\paragraph{Alternating methods for convex constraints} 

Due to its simplicity and effectiveness, alternating minimization has long been a popular optimization method, dating back to early work in the optimization literature (e.g. \citet{ortega1970iterative}), and has been widely studied under various assumptions (e.g.~\citet{auslender1976optimisation}, \citet{luo1993error}). 
For instance, assuming that the loss function $\loss$ is $\beta$-smooth and $\alpha$-strongly convex in each variable,
\citet{luo1993error} prove linear convergence for alternating minimization under convex constraints (and, in the case of more than two variables, for the analogous coordinate descent algorithm).

In some settings, the loss function $\loss$ may be more well-behaved with respect to one of the variables than the other, in terms of its smoothness and convexity properties.
\citet{beck2015convergence} studies alternating minimization for a convex loss $\loss(x,y)$ under convex constraints on $x$ and on $y$, 
proving that the gap in the loss function values, i.e.~the difference $\loss(x_t,y_t) -\loss(\xh,\yh)$, decays according to the rate $\O{\frac{\min\{\beta_x,\beta_y\}}{t}}$, where $\beta_x$ and $\beta_y$ represent the smoothness parameters of the loss $\loss$
with respect to the variables $x$ and $y$ respectively. Interestingly, this rate is controlled by the better of the two smoothness parameters---that is, the algorithm will converge
rapidly as long as {\em at least one} of the two smoothness parameters is bounded.

Our main results demonstrate an analogous phenomenon under an additional (restricted) strong convexity assumption---in this setting,
we find a {\em linear} convergence rate, with the convergence radius determined by $\min\left\{\frac{\beta_x}{\alpha_x},\frac{\beta_y}{\alpha_y}\right\}$,
where $\beta_x,\beta_y$ are smoothness parameters as before, while $\alpha_x,\alpha_y$ are the (restricted) strong convexity parameters with respect to $x$ and $y$, respectively.
That is, the linear convergence rate
depends on the better of the two condition numbers, while in \citet{beck2015convergence}'s result, without strong convexity,
the sublinear convergence rate depends on the better of the two smoothness parameters.
Thus, while a main focus of our work is to establish convergence results in a nonconvex setting, 
even in the convex setting our results reveal the interesting role of the two relative condition numbers (i.e.~for the $x$ and the $y$ variables)
in determining the overall convergence rate.

\section{Optimization over nonconvex constraints}\label{sec:nonconvex_constraint}

In this section, we briefly review the notion of local concavity coefficients introduced in \citet{barber2017gradient}, measuring concavity of the set at any given point. 

\subsection{Local concavity coefficients}

One main challenge of working over nonconvex regions is that, since $\Xset$ and $\Yset$ are potentially nonconvex sets, the standard first-order optimality conditions under convex setting do not apply. Specifically, fixing any $y\in\Yset$ and defining
\[x_y = \argmin \{\loss(x,y):x\in\Xset\},\]
the nonconvexity of $\Xset$ means that we cannot assume that $\inner{x - x_y}{\nabla_x\loss(x_y,y)}\geq 0$ for all
other $x\in\Xset$ (and same when we reverse the roles of $x$ and $y$). This makes the analysis of optimization problem with nonconvex constraints difficult, 
since the first-order optimality condition is crucial for understanding convergence behavior.

In order to overcome this obstacle, \citet{barber2017gradient} recently proposed the notion of {\em local concavity coefficients} for any nonconvex set,
related to the notion of prox-regular sets in the analysis literature (\citet{federer1959curvature,colombo2010prox}).
These concavity coefficients measure the extent to which the set deviates from convexity using four different properties of a convex set, and prove that these multiple definitions are all equivalent. Here we consider the {\em curvature condition}, which can be seen as a natural relaxation of the geometric characterization of a convex set:

\begin{definition}\label{def:curvature} {\bf (Curvature condition.)}
	Let $\Zset\subset \R^d$ be a closed subset, containing a point $z\in\Zset$. We say that  $\Zset$ satisfies
	the {\em curvature condition} with respect to a norm $\norm{\cdot}$ at the point $z$ with parameter $\gamma_z$ if, for all $z'\in\Zset$,
	\begin{equation*}\label{eqn:curvature}
	\limsup_{t\rightarrow 0} \frac{\min_{w\in\Zset}\norm{w - ((1-t)z + tz')}}{t} \leq \gamma_z\norm{z - z'}_2^2.\end{equation*}
\end{definition}
The norm $\norm{\cdot}$ used in the definition of the curvature condition is not necessarily the $\ell_2$ norm, and may be chosen
to suit the problem at hand---it is intended to reflect the natural structure arising in the setting of the problem.
For example, in high-dimensional setting, we may instead choose a more structured norm such as the $\ell_1$ norm ($\norm{\cdot}=\norm{\cdot}_1$) for a sparsity-inducing constraint or the nuclear norm (the sum of the singular
values of a matrix) for a low-rank constraint ($\norm{\cdot}=\nucnorm{\cdot}$).

Based on \defref{curvature}, we next define the local concavity coefficient $\gamma_z(\Zset)$ for any point $z\in\Zset$. We write $\pr{\Zset}$
to denote (possibly non-unique) projection to $\Zset$ with respect to the $\ell_2$ norm.
\begin{definition}\label{def:lcc} {\bf (Concavity coefficients.)}
	For a closed subset $\Zset\subset \R^d$, let $\mathcal{D}\subset \Zset$ be a set of {\em degenerate points},
	\[
	\mathcal{D} =\left\{z\in\Zset: \textnormal{$\pr{\Zset}$ is not continuous in any neighborhood of $z$} \right\}.
	\]
	The {\em local concavity coefficient} $\gamma_z(\Zset)$ at $z\in\Zset$, with respect to a norm $\norm{\cdot}$ and its dual norm $\norm{\cdot}^*$, is then given by
	\begin{equation*}\label{eqn:lcc}
	\gamma_z(\Zset) =\begin{cases}
	\infty, & z \in \mathcal{D},\\
	\min\{\gamma_z \in[0,\infty]: \textnormal{Condition \eqnref{curvature} holds at $z\in\Zset$ with $\gamma_z$} \}, & z \notin \mathcal{D}.
	\end{cases}
	\end{equation*}
	The {\em global concavity coefficient} is given by
	\[\gamma(\Zset) = \sup_{z\in\Zset}\gamma_z(\Zset).\]
\end{definition}

Note that any convex set $\Zset$ trivially satisfies $\gamma_{z}(\Zset)=0$ for all $z\in\Zset$. Moreover, it has been shown in~ \citep{barber2017gradient}  that
these coefficients are easy to compute or bound for  many nonconvex sets that are commonly used in high-dimensional statistical models---for instance,
for the rank-constrained set $\Zset=\{Z\in\R^{d_1\times d_2}:\rank(Z)\leq r\}$, the coefficients are given by $\gamma_z(\Zset) = \frac{1}{2\sigma_r(Z)}$, where $\sigma_r(Z)$
is the $r$th singular value of the matrix $Z$.
Intuitively, the curvature condition ensures that while $\Zset$ may have nonconvex boundaries in general, this nonconvexity must be fairly ``smooth'' wherever the coefficient $\gamma_z(\Zset)$ is small. 
For more details on local concavity coefficient and its application to nonconvex optimization, see \citet{barber2017gradient}. The connection to the notion of  prox-regular sets from the nonsmooth analysis literature is also discussed in {\citep[Section 2.3]{barber2017gradient}}.

One useful property of the local concavity coefficients $\gamma_z(\Zset)$ is that
they equivalently characterize the extent to which the usual first-order optimality conditions are violated when minimizing over the set $\Zset$:
\begin{lemma}[{\citep[Theorem 2]{barber2017gradient}}]\label{lem:firstorder} {\bf (First-order optimality.)}
	Let $\gamma_z(\Zset)$ be a local concavity coefficient, as in \defref{lcc}, with respect to a norm $\norm{\cdot}$ and its dual $\norm{\cdot}^*$. Then, for any differentiable function $f:\R^d\rightarrow \R$ such that $z$ is a local minimizer of $f$ over $\Zset$,
	\begin{equation*}
	\inner{z'-z}{\nabla f(z)}\geq - \gamma_z(\Zset)\norm{\nabla f(z)}^*\norm{z'-z}^2_2 \; \textnormal{ for all $z'\in\Zset$}.\end{equation*}
\end{lemma}
Comparing to a convex setting, where first-order optimality properties ensure that $\inner{z'-z}{\nabla f(z)}\geq 0$ whenever $z$ is a local minimizer of $f$ over $\Zset$,
we see that a small coefficient $\gamma_z(\Zset)$ ensures that $\Zset$ behaves ``almost'' like a convex set in this regard.

Returning to the alternating minimization setting, we first fix norms $\normx{\cdot}$ and $\normy{\cdot}$ for the $x$ and $y$ variables---for instance,
for a low-rank + sparse problem, we might choose $\normx{\cdot}$ and $\normy{\cdot}$ to be the nuclear norm and the $\ell_1$ norm, respectively.
To simplify our exposition, we will assume  that 
our structured norms $\normx{\cdot},\normy{\cdot}$ are scaled 
to satisfy $\normx{\cdot},\normy{\cdot}\geq\norm{\cdot}_2$, which is  the case for many of the structured norms that arise in various applications
(such as the $\ell_1$ norm and nuclear norm).

Let the local concavity coefficients $\gamma_x(\Xset)$ and $\gamma_y(\Yset)$ be defined with respect to these potentially different norms.
\lemref{firstorder} allows us to obtain approximate first-order optimality conditions for the steps of the alternating minimization
algorithm---for instance, letting $x_y$ be a local minimum of the problem $\min\{\loss(x,y):x\in\Xset \}$
(i.e.~the $x$ update step of alternating minimization), then for all $x \in\Xset$,
\begin{equation}\label{eqn:firstorder}
\inner{x - x_y}{\nabla_x \loss(x_y,y)} \geq -\gamma_{x_y}(\Xset)\normx{\nabla_x\loss(x_y,y)}^*\norm{x - x_y}_2^2,
\end{equation}
and similarly for $y$. These bounds provide a critical ingredient for our convergence analysis.

\section{Convergence analysis of alternating minimization}\label{sec:altmin}

We now turn to our convergence result on the alternating minimization method. Given the loss function $\loss(x,y)$ which is differentiable,  we consider an optimization problem,
\[\textnormal{Minimize }\loss(x,y)\textnormal{ over }x\in\Xset,y\in\Yset,\]
where the sets $\Xset\subset\R^{d_x}$ and $\Yset\subset\R^{d_y}$ represent the structural constraints on the variables $x$ and $y$ respectively. 

Let $(\xh,\yh)$ be the target of our optimization problem, which formally we require only to be a {\em local} minimizer of $\loss(x,y)$---this is because $\loss(x,y)$ may
potentially be highly nonconvex or degenerate in regions $(x,y)$ far from the origin, and we may even have $\lim \loss(x,y)=-\infty$ as $(x,y)$ tends to infinity in some direction.
If this is the case, then the steps of the alternating minimization algorithm could potentially diverge, and it may instead be necessary
to choose our update steps locally.

To formalize this, define new constraint sets $\Xset_0 = \Xset\cap \ball_2(x_0,\rho_x)$ and $\Yset_0 = \Yset\cap\ball_2(y_0,\rho_y)$,
where $(x_0,y_0)$ is our initialization point.
These neighborhoods of the original constraint sets $\Xset$ and $\Yset$ are assumed to be sufficiently large so as to contain the target point $(\xh,\yh)$
(in other words, our initialization point $(x_0,y_0)$ was chosen to be close to the target $(\xh,\yh)$),
but sufficiently small so that the loss function $\loss(x,y)$ is well-behaved over this small region $\Xset_0\times\Yset_0$.

We then define
\[(\xh,\yh) = \argmin\left\{\loss(x,y) : x\in\Xset_0,y\in\Yset_0 \right\},\]
and run the alternating minimization algorithm locally by iterating the steps
\begin{equation}\label{eqn:altmin}
\begin{cases} x_{t} = \argmin_{x\in\Xset_0}\loss(x,y_{t-1}),\\ y_{t}=\argmin_{y\in\Yset_0}\loss(x_{t},y).\end{cases}\end{equation}
For our intuition, we should interpret these radius constraints, i.e.~working in $\Xset_0$ and $\Yset_0$ rather than in $\Xset$ and $\Yset$, as 
a technicality for the theory, which we do not need to actually implement in practice. In particular,
for many settings, the alternating minimization steps are implemented with some kind of local search 
procedure, such as gradient descent in the $x$ or the $y$ 
variable, which will move towards a nearby local minimizer without enforcing a radius constraint.
In other settings, even the {\em global} minimizer for the $x$ or the $y$ variable (while the other variable is fixed),
stays within a small neighborhood, without enforcing a radius constraint. 
In other words, the radius constraint will generally not be active, 
and thus we can often ignore it in our implementation of the algorithm.
However, for the theoretical results obtained here, we require it in order to be able to handle a broader
range of problems.

The following lemma proves that, if the radii $\rho_x, \rho_y$ are chosen to be small, the curvature conditions (\defref{curvature}) of $\Xset$ and $\Yset$ are  inherited by $\Xset_0$ and $\Yset_0$:
\begin{lemma}\label{lem:LCC_rho}
	If $\rho_x< \frac{1}{2\max_{x\in\Xset_0}\gamma_x(\Xset)}$, then 
	$\gamma_x(\Xset_0)\leq \gamma_x(\Xset)$
	for all $x\in\Xset_0$, and in particular,
	\[\gamma(\Xset_0)\coloneqq  \sup_{x\in\Xset_0}\gamma_x(\Xset_0)  \leq \sup_{x\in\Xset_0}\gamma_x(\Xset).\]
	The analogous statement holds for $y$. 
\end{lemma}
The proof of this result is given in \appref{lem_LCC_rho}.

To see how this result will play a role in the convergence analysis for alternating minimization,
consider a single update step for the $x$ variable. Let $x_y = \argmin\{ \loss(x,y):x\in\Xset_0\}$. Then
\lemref{firstorder} proves the following bound (which we can compare to~\eqnref{firstorder}),
\begin{equation}\label{eqn:firstorder0}
\inner{x' - x_y}{\nabla_x \loss(x_y,y)} \geq -\gamma(\Xset_0)\normx{\nabla_x\loss(x_y,y)}^*\norm{x' - x_y}_2^2 \; \textnormal{ for all $x' \in\Xset_0$},
\end{equation}
while \lemref{LCC_rho} proves a useful bound for $\gamma(\Xset_0)$ as long as $\rho_x$ is sufficiently small (and similarly for $y$).


\subsection{Assumptions}\label{sec:assumptions}

Next we formally establish our assumptions on the loss function $\loss(x,y)$ as well as initialization condition. 

\paragraph{Loss function}
We first define some notation. 
Our convergence results will be derived in terms of the {\em first-order divergence}, a measure of distance to the optimal points $\xh$ and $\yh$
that is defined relative to the loss function:
\begin{align}
&\diver{x}{\xh}=\inner{x-\xh}{\nabla_x\loss(x,\yh) - \nabla_x\loss(\xh,\yh)}, \;\textnormal{ and } \label{eqn:diverx}\\ &\diver{y}{\yh}=\inner{y-\yh}{\nabla_y\loss(\xh,y)-\nabla_y\loss(\xh,\yh)}. \label{eqn:divery}\end{align}
This divergence  has been used also in \citet{loh2013regularized} to prove statistical errors of any local minimum in the sparse regression setting. 
Note that, if $\loss$ is nonconvex, then potentially $\diver{x}{\xh}$ or $\diver{y}{\yh}$ may be negative.
Abusing notation, we define the square root of the divergence as
\[\sdiver{x}{\xh} = \sqrt{\max\left\{0,\diver{x}{\xh}\right\}} \;\; \textnormal{ and } \;\;
\sdiver{y}{\yh} = \sqrt{\max\left\{0,\diver{y}{\yh}\right\}},\]
to accommodate the case where the divergences may be negative.

Throughout we will write $\errx,\erry \geq 0$ to indicate vanishing error terms that allow a small amount of slack
in the convexity and smoothness conditions. In the high-dimensional statistics literature, these terms often represent the ``statistical error''---meaning, 
if the global minimizer $\xh$ approximates some ``true'' parameter $x^\star$ only up to an error level of $\errx$, then as soon as our iterative algorithm
reaches a solution $x_t$ within distance $\sim \errx$ of $\xh$, we are already optimal (up to a constant) in terms of estimating the underlying parameters $x^\star$.
While our work in this paper is not based in a concrete statistical model,  we will still refer to $\errx,\erry$ as the statistical error terms,
as this is often the case for many of the applications of our result.

We now state our assumptions on the loss $\loss(x,y)$. As mentioned earlier, our optimization method works locally in neighborhoods of the initialization point $(x_0,y_0)$. Consequently, it is sufficient for us to require the assumptions on $\loss(x,y)$ to hold only locally in the regions $\Xset_0$ and $\Yset_0$. 

First, since $(x,y)$ are being optimized jointly, we need to ensure that these two variables are identifiable, and require a joint restricted strong convexity (RSC) condition at the target point $(\xh,\yh)$:
\begin{assumption}\label{assump:RSC_joint} {\bf (Joint restricted strong convexity (RSC).)}
	For all $x\in\Xset_0$ and all $y\in\Yset_0$,
	\begin{equation}\label{eqn:RSC_joint}\inner{\left(\begin{array}{c}x - \xh\\ y - \yh \end{array}\right)}
	{\nabla\loss(x,y)- \nabla\loss(\xh,\yh)}
	\geq \alpha_x\norm{x-\xh}^2_2  + \alpha_y\norm{y-\yh}_2^2 - \alpha_x\errx^2 - \alpha_y\erry^2.
	\end{equation}
\end{assumption}
Note that we require joint RSC to hold only at the target $(\xh,\yh)$. In other regions of $\Xset \times\Yset$, we may not have joint convexity if the variables $x$ and $y$ are not identifiable from each other
in general (for instance, this arises in low-rank + sparse decomposition problems). 

Next, we assume that, marginally in $x$ and in $y$, the loss function satisfies the restricted smoothness (RSM) property near the optimal point $(\xh,\yh)$:
\begin{assumption}\label{assump:RSM} {\bf (Restricted smoothness (RSM).)}
	For all $x\in\Xset_0$ and all $y\in\Yset_0$,
	\begin{equation}\label{eqn:RSM_marg}\diver{x}{\xh}\leq 
	\beta_x\norm{x-\xh}_2^2 +\alpha_x\errx^2\textnormal{\quad and \quad}
	\diver{y}{\yh}\leq 
	\beta_y\norm{y-\yh}_2^2 +\alpha_y\erry^2.
	\end{equation}
\end{assumption}
Comparing to the restricted strong convexity assumption, we see that we need to choose constants 
$\alpha_x\leq\beta_x$ and $\alpha_y\leq\beta_y$. 

We also mention that we allow one of the smoothness parameters to be much larger relative to the other, i.e.~$\beta_x \gg \beta_y$ or vice versa, unlike methods that work jointly in the combined $(x,y)$ variable (e.g.~gradient descent on this single combined
variable), whose performance is closely tied to the smoothness of the total problem, $\max\{\beta_x,\beta_y\}$.
We will discuss this distinction in more detail in \secref{conditioning}.


Finally, we require a ``cross-product'' condition (explained below):
\begin{assumption}\label{assump:crossprod} {\bf (Cross-product bound.)}
	For all $x\in\Xset_0$ and all $y\in\Yset_0$,
	\begin{multline*}\left|\inner{x-\xh}{\nabla_x\loss(x,y) - \nabla_x\loss(x,\yh)} - \inner{y-\yh}{\nabla_y \loss(x,y) - \nabla_y \loss(\xh,y)}\right| \\
	\leq 
	\frac{1}{2}\mu_x\norm{x-\xh}^2_2 +\frac{1}{2} \mu_y\norm{y-\yh}_2^2 +\alpha_x\errx^2 + \alpha_y  \erry^2,\end{multline*}
	where $0\leq \mu_x\leq \alpha_x$ and $0\leq \mu_y\leq \alpha_y$.
\end{assumption}
To understand this assumption, suppose that $\loss$ is twice differentiable. In this case, applying Taylor's theorem to rewrite the above
expression in terms of $\nabla^2\loss$, we find that \assumpref{crossprod} holds with
\[\mu_x = \mu_y = \sup_{\substack{x\in\Xset_0;y\in\Yset_0\\t,t'\in[0,1]}}2\opnorm{\nabla^2_{xy}\loss(x,ty+(1-t)\yh) - \nabla^2_{xy}\loss(t'x+(1-t')\xh,y)},\]
where the norm $\opnorm{\cdot}$ is the matrix operator norm (the largest singular value). Since $\Xset_0$ and $\Yset_0$ are bounded via the radii $\rho_x,\rho_y$, then, this condition is satisfied
whenever $\nabla^2_{xy}$ is Lipschitz. As a special case, if $\loss(x,y)$ is quadratic, then we can trivially take $\mu_x=\mu_y=0$
since $\nabla^2_{xy}$ is constant.
%

\paragraph{Initialization}

As our theoretical results mainly concern the local behavior of the alternating minimization method, the initialization scheme is crucial to ensure the success of the procedure. 
Our results require the following  initialization condition:

\begin{assumption}\label{assump:initialization} {\bf (Initialization condition.)} Let $\gamma(\Xset_0)$ and $\gamma(\Yset_0)$ be the concavity coefficients of $\Xset_0$ and $\Yset_0$ with respect to norms $\normx{\cdot}$ and $\normy{\cdot}$ (and its duals $\normx{\cdot}^*$ and $\normy{\cdot}^*$) respectively. Then
	\[2\gamma(\Xset_0)\cdot  \left(\normx{\nabla_x\loss(\xh,\yh)}^* + \max_{y\in\Yset_0}\normx{\nabla_x\loss(x_y,y)}^* \right)  \leq \alpha_x - \mu_x,\]
	and
	\[2\gamma(\Yset_0)\cdot\left(\normy{\nabla_y\loss(\xh,\yh)}^* + \max_{x\in\Xset_0}\normy{\nabla_y\loss(y,y_x)}^*\right) \leq \alpha_y - \mu_y,\]
	where $x_y=\argmin\{\loss(x,y):x\in\Xset_0 \}$ and $y_x$ is defined similarly.
\end{assumption}

Recall that \lemref{LCC_rho} provides easy bounds on $\gamma(\Xset_0)$ and $\gamma(\Yset_0)$, as long as the radii
$\rho_x,\rho_y$ are chosen to be sufficiently small; furthermore, if $\Xset$ is convex, then $\gamma(\Xset_0)=0$ and so the first
bound holds trivially, and similarly for the second bound if $\Yset$ is convex.
In the nonconvex setting where $\gamma(\Xset_0)$ and/or $\gamma(\Yset_0)$ are nonzero, 
see \citep[Section 3.4]{barber2017gradient} for a discussion of the necessity of this type of initialization condition
for the related problem of gradient descent in a single variable; we believe that this type of condition is 
necessary for alternating minimization as well, in the absence of additional assumptions.

\subsection{Convergence guarantee}

Now we show that by alternating optimization over $x$ and over $y$, we obtain fast convergence to the target $(\xh,\yh)$ up to the level of a small statistical error term. We prove convergence by working with the first-order divergence defined in~\eqnref{diverx},~\eqnref{divery} above.

While the divergence may take negative values in general, according to \assumpref{RSC_joint}, it will be always nonnegative in the regions $\Xset_0$ and $\Yset_0$, up to the statistical error. 
The following result then provides guarantee on convergence of alternating minimization \eqnref{altmin} as measured in the square-root divergences $\sdiver{x}{\xh}$ and $\sdiver{y}{\yh}$:

\begin{theorem}\label{thm:altmin_contract}
	Suppose that 
	\assumpsref{RSC_joint},~\assumpssref{RSM},~\assumpssref{crossprod}, 
	and~\assumpssref{initialization} hold.
	Then the iterations of the alternating minimization algorithm~\eqnref{altmin} satisfy the recursive bounds
	\begin{align}
	&\label{eqn:thm1_x}\sdiver{x_t}{\xh} \leq \sqrt{1-\frac{\alpha_y}{2\beta_y} } \cdot\sdiver{y_{t-1}}{\yh}+ \sqrt{3(\alpha_x\errx^2+\alpha_y\erry^2)},\textnormal{ 
		and}\\
	&\label{eqn:thm1_y}\sdiver{y_t}{\yh} \leq \sqrt{1- \frac{\alpha_x}{2\beta_x} } \cdot\sdiver{x_t}{\xh}+ \sqrt{3(\alpha_x\errx^2+\alpha_y\erry^2)},\end{align}
	for all $t\geq 1$. In particular, this implies a linear rate of convergence:
	\begin{equation}\label{eqn:converge_L2}
	\norm{(x_t,y_t) - (\xh,\yh)}_2 \leq\left(\sqrt{1 - \frac{\alpha_x}{2\beta_x}}\cdot \sqrt{1 - \frac{\alpha_y}{2\beta_y}}\right)^t \cdot \frac{\sqrt{6\beta_y}\rho_y}{\sqrt{\min\{\alpha_x,\alpha_y\}}} + C\cdot\max\{\errx,\erry\}\end{equation}
	for all $t\geq 1$, where
	\[ C =  \frac{18}{1-\sqrt{1 - \frac{\alpha_x}{2\beta_x}}\cdot\sqrt{1 - \frac{\alpha_y}{2\beta_y}}}\cdot \sqrt{\frac{\max\{\alpha_x,\alpha_y\}}{\min\{\alpha_x,\alpha_y\}}}.\]
\end{theorem}
This theorem is proved in \secref{proofs_thms}.

Before proceeding, we remark that the order of the updates---that is, after initializing at time $t=0$ with points $x_0,y_0$,
at time $t=1$ we then update first $x$ and then $y$---is arbitrary. In particular,
the term $\sqrt{\beta_y} \rho_y$ appearing in the numerator of~\eqnref{converge_L2},
can of course be replaced instead by $\sqrt{\beta_x}\rho_x$ if we switch the order of the updates. 
This suggests that it  may be best to first update the variable with {\em poorer} smoothness parameter---that is, 
if the $y$ variable is more well-conditioned, at our first step we should fix $y$ and update $x$.

\subsubsection{Dependence on condition number}\label{sec:conditioning}
Examining the bound~\eqnref{converge_L2} for the convergence rate in the $\ell_2$ norm,
we see that the convergence rate is dominated by the radius
\[\sqrt{1 - \frac{\alpha_x}{2\beta_x}}\cdot \sqrt{1 - \frac{\alpha_y}{2\beta_y}}\]
(here we ignore the negligible statistical error term $C\cdot \max\{\errx,\erry\}$).
We now discuss the implications of this result, in terms of its dependence on the convexity and smoothness parameters,
$\alpha_x,\alpha_y$ and $\beta_x,\beta_y$.
To help us discuss the conditioning of this problem, we define the two marginal condition numbers 
of the loss function with respect to the $x$ and the $y$ variables,
\[\kappa_x(\loss) = \frac{\beta_x}{\alpha_x}\textnormal{ and }\kappa_y(\loss) = \frac{\beta_y}{\alpha_y},\]
and the joint condition number
\[\kappa(\loss) = \frac{\max\{\beta_x,\beta_y\}}{\min\{\alpha_x,\alpha_y\}} \geq \max\{\kappa_x(\loss),\kappa_y(\loss)\},\]
which, up to constant factors, gives the condition number of the loss function $\loss$ as a function of the joint variable $(x,y)$.

In~\eqnref{converge_L2}, we see that our convergence radius is strictly smaller than $1$, as long
as {\em either} of the two marginal condition numbers
is bounded from above, that is, if $\min\{\kappa_x(\loss),\kappa_y(\loss)\}$ is bounded.
On the other hand, if we consider optimization algorithms that work with the combined joint variable $(x,y)$,
the performance of such algorithms typically relies heavily on the joint
condition number $\kappa(\loss)\geq \max\{\kappa_x(\loss),\kappa_y(\loss)\}$.
For example, if $\loss$ is $\alpha$-strongly convex and $\beta$-smooth in the joint variable $(x,y)$, standard results (see e.g. \citet{bubeck2015convex}) prove that gradient descent in $(x,y)$ yields
\[\norm{(x_t,y_t) - (\xh,\yh)}_2\leq \big(\sqrt{1-\alpha/\beta}\big)^t \norm{(x_0,y_0) - (\xh,\yh)}_2.\]
Comparing to our notation, it can be shown that $\alpha\leq\min\{\alpha_x,\alpha_y\}$ and $\beta\geq\max\{\beta_x,\beta_y\}$, 
and so the radius of covergence for (joint) gradient descent is controlled by the joint condition number, $\kappa(\loss)\geq \max\{\kappa_x(\loss),\kappa_y(\loss)\}$.

Therefore, in settings where one of the two---$\kappa_x(\loss)$ or $\kappa_y(\loss)$---is much larger than the other,
we may expect that (joint) gradient descent, or other non-alternating algorithms, might perform poorly, 
while alternating minimization will continue to perform well, since its linear convergence rate depends
only on the best of the two condition numbers, i.e.~on $\min\{\kappa_x(\loss),\kappa_y(\loss)\}$.
(As discussed earlier in \secref{relatedwork},
\citet{beck2015convergence} find an analogous result without strong convexity assumptions, demonstrating that the sublinear rate of convergence for alternating minimization method is driven by minimum of the two smoothness parameters, i.e. $\min\{\beta_x,\beta_y\}$.)

We will explore this phenomenon empirically when we present numerical experiments with simulated data (see \secref{empirical} below).

\section{Inexact alternating minimization}\label{sec:inexact}
In some settings, it may be impractical to solve the alternating minimization steps exactly, i.e.~when $\loss(x,y)$
is difficult to minimize even as a function of only $x$ or only $y$. 
In these cases, we may want to solve each step of the alternating minimization algorithm inexactly.
We first state a general result for this inexact setting, then discuss the specific strategy
of taking the approximate steps via alternating gradient descent.

To study the convergence behavior of alternating minimization where the steps are computed only approximately,
we formulate an inexact algorithm where, at each step, we choose
$x_{t}$ and $y_{t}$ to be within some tolerance parameters $\delta^x_t,\delta^y_t$ of the
exact alternating minimization steps at that time: for all $t\geq 1$, 
\begin{equation}\label{eqn:altmin_inexact}
\begin{cases} x^{\textnormal{exact}}_{t} = \argmin_{x\in\Xset_0}\loss(x,y_{t-1}),\quad x_{t}\in\Xset_0 \cap \ball_2(x_{t}^{\textnormal{exact}},\delta^x_{t}),\\ y^{\textnormal{exact}}_{t}=\argmin_{y\in\Yset_0}\loss(x_{t},y), \quad y_{t}\in\Xset_0 \cap \ball_2(y_{t}^{\textnormal{exact}},\delta^y_{t}).\end{cases}\end{equation}
Here $x_t$ and $y_t$ can be chosen arbitrarily (or even adversarially) as long as they are within the required distance
of the true solutions $x^{\textnormal{exact}}_t$ and $y^{\textnormal{exact}}_t$.

In order to establish the convergence of the inexact alternating minimization algorithm \eqnref{altmin_inexact}, we require
an additional assumption: 
\begin{assumption}\label{assump:triangle} {\bf (Relaxed triangle inequality.)}
	For all $x,x'\in\Xset_0$,
	\[\sdiver{x}{\xh}\leq \sdiver{x'}{\xh} + \sqrt{\beta_x}  \norm{x-x'}_2 +\sqrt{\alpha_x}  \errx,\]
	and for all $y,y'\in\Yset_0$,
	\[\sdiver{y}{\yh}\leq \sdiver{y'}{\yh} + \sqrt{\beta_y} \norm{y-y'}_2 + \sqrt{\alpha_y} \erry,\]
\end{assumption}
It can be shown that a stronger form of the restricted smoothness condition (\assumpref{RSM}) implies this type of 
relaxed triangle inequality, but for simplicity we state it as an assumption.

The following theorem states that the inexact alternating minimization inherits fast convergence of the alternating minimization steps to the target $(\xh,\yh)$,
under the same assumptions as the original result \thmref{altmin_contract}, along with the relaxed triangle inequality (\assumpref{triangle}).
\begin{theorem}\label{thm:altmin_inexact}
	Suppose that 
	\assumpsref{RSC_joint},~\assumpssref{RSM},~\assumpssref{crossprod},~\assumpssref{initialization}, and~\assumpssref{triangle} hold. Then,
	the steps of the inexact alternating minimization algorithm satisfy
	\begin{align}
	&\label{eqn:thm2_x}
	\sdiver{x_t}{\xh}\leq \sqrt{1 - \frac{\alpha_y}{2\beta_y}}\cdot \sdiver{y_{t-1}}{\yh} +  \sqrt{\beta_x} \delta_t^x+ \sqrt{8(\alpha_x\errx^2+\alpha_y\erry^2)} \textnormal{\quad and}\\
	&\label{eqn:thm2_y}
	\sdiver{y_t}{\yh}\leq \sqrt{1 - \frac{\alpha_x}{2\beta_x}}\cdot \sdiver{x_t}{\xh} +  \sqrt{\beta_y} \delta_t^y + \sqrt{8(\alpha_x\errx^2+\alpha_y\erry^2)},
	\end{align}
	for all $t\geq 1$.
\end{theorem}
This theorem is proved in \secref{proofs_thms}.

Of course, in order for this result to be meaningful, the slack terms $\delta_t^x,\delta_t^y$ need to be sufficiently small,
so that the errors $\sdiver{x_t}{\xh}$ and $\sdiver{y_t}{\yh}$ are able to converge to zero (or, at least, to the level of the statistical error
terms $\errx,\erry$).
As a special case, consider the setting where the slack terms $\delta_t^x,\delta_t^y$ decrease
as the solution converges, via the rule
\begin{equation}\label{eqn:delta_rule_recursive}
\delta_t^x\leq c_x\norm{x_{t-1}-x_t^{\textnormal{exact}}}_2+ C_x\errx, \quad \delta_t^y\leq c_y\norm{y_{t-1}-y_t^{\textnormal{exact}}}_2+ C_y\erry,
\end{equation}
for some sufficiently small $c_x,c_y\geq 0$ and for some $C_x,C_y<\infty$. 

In fact, we will see momentarily that this recursive bound arises
naturally when the approximate iterative solutions $x_t$ and $y_t$ are obtained via alternating gradient descent.
First, however, we prove that the recursive rule~\eqnref{delta_rule_recursive} is sufficient
to ensure linear convergence as long as the constants $c_x,c_y$ are sufficiently small.

\begin{lemma}\label{lem:altmin_inexact}
	Suppose that, for all $t\geq 1$, the slack terms $\delta_t^x,\delta_t^y$ satisfy~\eqnref{delta_rule_recursive}.
	Then, under the assumptions of \thmref{altmin_inexact}, 
	if
	\begin{equation}\label{eqn:radius_inexact}
	r\coloneqq \left( \sqrt{1 - \frac{\alpha_x}{2\beta_x}}  + 3c_y\sqrt{\frac{\beta_y}{\alpha_y}}\right) \cdot\left(\sqrt{1 - \frac{\alpha_y}{2\beta_y}}  + 3c_x\sqrt{\frac{\beta_x}{\alpha_x}}\right) < 1,\end{equation}
	then the iterations of the inexact alternating minimization algorithm~\eqnref{altmin_inexact} satisfy
	\[
	\norm{(x_t,y_t) - (\xh,\yh)}_2 \leq r^t \cdot \frac{\sqrt{6(\alpha_x\rho_x^2 + \beta_y \rho_y^2)}}{\sqrt{\min\{\alpha_x,\alpha_y\}}} + C\cdot \max\{\errx,\erry\}\]
	for all $t\geq 1$, where 
	\[C = \frac{39}{1-r}\cdot \sqrt{\frac{\alpha_x + \alpha_y + C_x^2\beta_x + C_y^2\beta_y}{\min\{\alpha_x,\alpha_y\}}}.\]
\end{lemma}
The proof of this lemma is given in \appref{lem_altmin_inexact}.

We should interpret this lemma as covering two distinct scenarios:
\begin{itemize}
	\item[$\bullet$] First, if the loss is well-conditioned in both the $x$ and the $y$ variables---that is, both $\frac{\beta_x}{\alpha_x}$ and $\frac{\beta_y}{\alpha_y}$ are bounded---then
	we can afford inexact update steps for both variables, allowing $c_x,c_y$ to both be small positive constants while still obtaining linear convergence.\smallskip
	\item[$\bullet$] Alternately, if the loss is well-conditioned in one variable only---without loss of generality, if $\frac{\beta_x}{\alpha_x}$ is large (or even $\beta_x=\infty$)
	while $\frac{\beta_y}{\alpha_y}$ is bounded---then we can allow the $y$ variable update to be performed inexactly,
	while the $x$ variable should be updated with the exact alternating minimization step (that is, $c_x=C_x=0$, i.e.~$\delta^t_x=0$ at each update iteration $t$).
	In this case, we can still obtain a linear convergence rate.

\end{itemize} 

In this second setting, to achieve linear convergence of an inexact algorithm, the exact update of the $x$ variable is critical. Analogous results have also appeared in the more general case of $n\geq 2$ variables, where it is known that without strong convexity assumptions, block-coordinate descent methods (e.g. \citet{nesterov2012efficiency}), which, at every iteration, performs a gradient descent on each variable (block), converges with the rate depending on the sum of the smoothness parameters of the $n$ variables; therefore, even for $n=2$ case, one cannot hope for fast convergence if one variable has poor smoothness, unless it is optimized exactly. On the other hand, as long as  exact minimization is performed on the least smooth variable, the convergence rate of block-coordinate descent  can be shown to scale independently of the least smoothness parameter (\citet{diakonikolas2018alternating}). 

While these types of results require that updating one variable is easier than the other, we will see examples of this setting in the simulations. More examples can  be found in e.g.~\citet{beck2015convergence,diakonikolas2018alternating,jain2015alternating}.

\subsection{Alternating gradient descent}\label{sec:altgrad}

If the alternating minimization update for $x$ is approximated via gradient descent in $x$ (and
same for $y$), then we may expect the errors in each step to scale linearly as in~\eqnref{delta_rule_recursive}.
We now give details for this claim,  relying on earlier work \citep{barber2017gradient}, which we summarize here:

To run the alternating descent algorithm, we first
initialize at $x_0\in\Xset$, $y_0\in\Yset$, then, at each iteration $t=1,2,\dots$,
\begin{enumerate}
	\item Perform $m_x$ many gradient descent steps for $x$:
	\begin{equation}\label{eqn:grad_descent_x}\!\!\!\!\!\begin{cases}
	\textnormal{Set $x_{t;0} = x_{t-1}$;}\\
	\textnormal{For $m=1,\dots,m_x$, set $x_{t;m} = \Pr{\Xset_0}{x_{t;m-1} - \eta_x\nabla_x \loss(x_{t;m-1},y_{t-1})}$;}\\
	\textnormal{Set $x_t = x_{t;m_x}$.}\end{cases}\end{equation}
	\item Perform $m_y$ many gradient descent steps for $y$:
	\begin{equation}\label{eqn:grad_descent_y}\!\!\!\!\!\begin{cases}
	\textnormal{Set $y_{t;0} = y_{t-1}$;}\\
	\textnormal{For $m=1,\dots,m_y$, set $y_{t;m} = \Pr{\Yset_0}{y_{t;m-1} - \eta_y\nabla_y \loss(x_t,y_{t;m-1})}$;}\\
	\textnormal{Set $y_t = y_{t;m_y}$.}\end{cases}\end{equation}
\end{enumerate}
Here $\pr{\Xset_0}$ and $\pr{\Yset_0}$ denote projection, with respect to the $\ell_2$ norm, to the sets $\Xset_0$ and $\Yset_0$. (Of course,
this projection may not be unique in the presence of nonconvexity---if this is the case, we can take any one of the closest points.)

Often, the number of steps in each ``inner loop'', namely $m_x$ and $m_y$, can be taken to be a small constant (or even $1$)
to obtain good empirical performance (note, in the case of $m_x=m_y=1$, the alternating descent algorithm corresponds to cycle block-coordinate descent with $n=2$ blocks). In our theoretical analysis,
under restricted strong convexity and restricted smoothness, constant values for $m_x$ and $m_y$ suffice to guarantee convergence.

\paragraph{Convergence results}

Gradient descent (and its variants) now serves as a popular tool in large-scale optimization problem, due to its scalablity to the high-dimensional setting. 
For convex constraints, the convergence behavior of the gradient descent has been well established (e.g. \citet{nesterov2013introductory}, \citet{nesterov2007gradient}, \citet{bubeck2015convex}) and also generalized to the high-dimensional setting (e.g. \citet{agarwal2010fast}). In contrast to the convex setting, gradient descent under nonconvex
constraints is more challenging to analyze theoretically, although it often performs well empirically.

In our previous work \citep{barber2017gradient}, we have shown that, under the framework of local concavity measure (\defref{lcc}), gradient descent converges rapidly to the optimal as long as it is initialized near the target. Here we follow the same setup, but slightly modify the assumptions to better match the setting of the alternating minimization problem.
We assume that:
\begin{itemize}
	\item[$\bullet$] The ``marginal'' loss functions $x\mapsto \loss(x,y)$ and $y\mapsto\loss(x,y)$ satisfy restricted strong convexity (RSC) uniformly over $\Xset_0$ and over $\Yset_0$,
	respectively. In other words, for all $x,x'\in\Xset_0$ and all $y\in\Yset_0$,
	\begin{equation}\label{eqn:RSC_obj} \loss(x,y) - \loss(x',y) - \inner{\nabla_x\loss(x',y)}{x-x'} \geq \frac{\alpha_x}{2}\norm{x-x'}_2^2 - \frac{\alpha_x}{2}\cdot\errx^2 , \end{equation}
	and analogously with the roles of $x$ and $y$ reversed.\smallskip
	
	\item[$\bullet$] The ``marginal'' functions $x\mapsto \loss(x,y)$ and $y\mapsto\loss(x,y)$ satisfy restricted smoothness (RSM) in $\Xset_0$ and $\Yset_0$, that is,
	\begin{equation}\label{eqn:RSM_obj}
	\loss(x,y) - \loss(x',y) - \inner{\nabla_x\loss(x',y)}{x-x'} \leq \frac{\beta_x}{2}\norm{x-x'}_2^2 + \frac{\alpha_x}{2}\cdot\errx^2,
	\end{equation}
	and analogously with the roles of $x$ and $y$ reversed.\smallskip
	
	\item[$\bullet$] There exist parameters $\phi_x, \phi_y \geq 1$ such that for any $z\in\R^{d_x}$, 
	\begin{equation}\label{eqn:compat}\normx{z-\pr{\Xset_0}(z)}^* \leq \phi_x \min_{x'\in\Xset_0}\normx{z-x'}^*,\end{equation}
	and analogously with the roles of $x$ and $y$ reversed (note that this condition holds trivially with $\phi_x=1$ if $\normx{\cdot}=\norm{\cdot}_2$,
	but $\phi_x=1$ often suffices even for other norms). \smallskip
	
	\item[$\bullet$] For all $y\in\Yset_0$,
	\begin{equation}\label{eqn:initialization2}2\phi_x\cdot  \max_{x,x'\in\Xset_0}\gamma_x(\Xset) \normx{\nabla_x\loss(x',y)}^* \leq (1-a_x)\cdot\alpha_x\end{equation}
	for some constant $a_x>0$, and analogously with the roles of $x$ and $y$ reversed with some constant $a_y>0$.
	
\end{itemize}

We remark that the assumptions~\eqnref{RSC_obj} and~\eqnref{RSM_obj} are made in terms of objective function, as opposed to the gradient forms
in the corresponding conditions from earlier, given in~\eqnref{RSC_joint} for RSC and~\eqnref{RSM_marg} for RSM. However, in 
general, the two forms of these conditions are roughly interchangeable with each other. 
For detailed discussion on the rest of the assumptions, see \citet{barber2017gradient} and references therein.

With these assumptions in place, we can bound the tolerance parameters $\delta_t^x, \delta_t^y$ appearing
in the inexact alternating minimization algorithm~\eqnref{altmin_inexact}, when the inexact steps are computed via gradient descent.
\begin{lemma}[{\citep[Theorem 3]{barber2017gradient}}]\label{lem:grad_descent}
	Suppose that  conditions~\eqnref{RSC_obj}, \eqnref{RSM_obj}, \eqnref{compat} and \eqnref{initialization2} hold. Then the output of 
	$m_x$ many gradient descent steps on the $x$ variable, given in~\eqnref{grad_descent_x}, satisfies
	\[\norm{x_t- x_{t}^{\textnormal{exact}}}_2^2 \leq \left(1 - a_x\frac{2\alpha_x}{\alpha_x + \beta_x}\right)^{m_x} \cdot \norm{x_{t-1} - x_t^{\textnormal{exact}}}_2^2 + \frac{1.5}{a_x} \cdot \errx^2. \]
	The analogous statement holds with the roles of $x$ and $y$ reversed.
\end{lemma}

Examining the conditions~\eqnref{delta_rule_recursive} and~\eqnref{radius_inexact} on the allowed size of the slack terms $\delta_t^x$ 
and $\delta_t^y$, we can see that taking $c_x,c_y = \O{\sqrt{\alpha_x/\beta_x}\cdot\sqrt{\alpha_y/\beta_y }}$ 
is sufficient to ensure that the condition~\eqnref{radius_inexact} will hold.
This yields the following corollary, which we state informally to avoid overly complicated constants:
\begin{corollary}\label{cor:approx_via_grad_descent}
	Under the assumptions of \lemsref{altmin_inexact} and \lemssref{grad_descent}, for some radius ${\textnormal{Rad}=\O{\sqrt{1-\frac{\alpha_x}{\beta_x}} \cdot \sqrt{1-\frac{\alpha_y}{\beta_y}}}<1}$,
	\[\norm{(x_t,y_t)-(\xh,\yh)}_2 \leq \O{\textnormal{Rad}^t \cdot \max\{\rho_x,\rho_y\} + \max\{\errx,\erry\}}\]
	for all $t\geq 1$ as long as
	\[\begin{cases}
\textnormal{Either the $x$ update  is  exact, or is approximated via $m_x = \O{\frac{\beta_x}{\alpha_x} \max\left\{ \log\left(\frac{\beta_x}{\alpha_x}\right),\log\left(\frac{\beta_y}{\alpha_y}\right) \right\}}$ many steps of gradient}\\\textnormal{  descent in $x$;}\\
	\textnormal{Either the $y$ update  is  exact, or is approximated via $m_y = \O{\frac{\beta_y}{\alpha_y} \max\left\{ \log\left(\frac{\beta_x}{\alpha_x}\right),\log\left(\frac{\beta_y}{\alpha_y}\right) \right\}}$ 
		many steps of gradient}\\\textnormal{ descent in $y$.}\end{cases}\]
\end{corollary}

To summarize, our results show that the convergence radius of the outer loops of alternating descent depends on the smaller of the two condition numbers, $\min\{\kappa_x(\loss),\kappa_y(\loss) \}$, while the number of steps in each inner loop (if performed inexactly) is mainly dependent on $\kappa_x$ for the $x$ variable and on $\kappa_y$ for the $y$ variable respectively (up to a logarithmic factor $\log \max\{\kappa_x(\loss),\kappa_y(\loss) \}$). Therefore, in the setting where both the condition numbers $\kappa_x$ and $\kappa_y$ are well bounded, the alternating descent algorithm still attains linear convergence rate, though the larger of the two condition numbers also becomes crucial in the total convergence time. Further, as long as one of the variables is updated exactly, the convergence rate remains unchanged even when its condition number is much larger relative to the other or even unbounded, i.e. $\kappa_x\gg \kappa_y$ or $\kappa_x=\infty$.

\section{Examples}\label{sec:examples}

This section highlights applications of our general theory to two classes of low-rank estimation problems, low-rank + sparse matrix decomposition and reduced rank multitask regression. In each of the settings, we 
present a target optimization problem and verify the assumptions of our main results, under the suitable choices of the radii $\rho_x,\rho_y$. All results in this section are proved in \appref{examples_proofs}.

\subsection{Matrix decomposition}\label{sec:matrix_decomposition}

In robust principal component analysis (RPCA) and in Gaussian factor models, it is common to assume that the matrix of interest is formed by a sum of low-rank and sparse components. Given the data generated through a matrix $\Xs + \Ys$, where $\Xs$ is low-rank and $\Ys$ is sparse, the task of matrix decomposition is to recover both the low-rank and sparse components simultaneously. 

While this problem is generally known to be ill-posed, certain incoherence conditions on the low-rank component have been shown to guarantee exact or approximate recovery of both low-rank and sparse components (e.g. \citet{candes2011robust,chandrasekaran2011rank}). Following \citet{negahban2012restricted}, we rely on the condition called {\em spikiness condition}, which restricts the class of low-rank matrices by imposing the constraint $\norm{\Xs}_\infty \leq \frac{\spike}{d}$.

Suppose that the underlying low-rank matrix $\Xs$ is positive semidefinite with rank $r$ and the sparse matrix $\Ys$ is symmetric with at most $s$ nonzero entries per each row. This gives rise to the following constraint sets:\footnote{For our analysis, we assume that $\norm{\Ys}_1$ is known {\em exactly}; this is a common 
	assumption for the constrained problem, e.g.~see \citet{amelunxen2014living}. On the other hand, $\norm{\Xs}_\infty$ needs only to be bounded by some known value $\spike/d$; we do not
	need to know it exactly.}
\begin{equation}\label{eqn:XY_lowrank_sparse}
\Xset = \left\{X\in\mathbb{S}_{+}^{d\times d}:\rank(X) \leq r, \; \norm{X}_\infty \leq \frac{\spike}{d}\right\}, \ 
\Yset = \left\{Y\in\mathbb{S}^{d\times d}:\norm{Y}_1\leq \norm{\Ys}_1\right\},
\end{equation}
where $\mathbb{S}$ and $\mathbb{S}_+$ denote the sets of symmetric or positive semidefinite $d\times d$ matrices, respectively. The $\ell_\infty$-ball constraint in the set $\Xset$ makes sure that the low-rank update by the algorithm is at most $\spike$-spiky at every iteration. To compute the local concavity coefficient of the set $\Xset$, the following lemma provides the upper bound on the concavity coefficient $\gamma_X(\Xset)$:
\begin{lemma}\label{lem:matrix_lcc}
	For the constraint set $\Xset=\{X\in\mathbb{S}_+^{d\times d}: \rank(X)\leq r, \norm{X}_\infty \leq \frac{\spike}{d}\}$, we have $\gamma_X(\Xset)\leq \frac{5}{4\sigma_r(X)}$ with respect to the nuclear norm $\norm{\cdot}_{X}=\nucnorm{\cdot}$.
\end{lemma}

Note that while the set of rank-constrained matrices without spikiness constraint has the local concavity coefficient $\gamma_X(\Xset) = \frac{1}{2\sigma_r(X)}$ (see \citep[Lemma 7]{barber2017gradient}), the lemma above shows that the coefficient for $\Xset$ can be upper bounded with a larger constant factor.

For a given loss function $\loss(X,Y)$,  we can then recover the underlying matrices $(\Xs,\Ys)$ by solving the constrained optimization problem 
\[(\Xh,\Yh)=\argmin\{ \loss(X,Y):(X,Y)\in\Xset\times\Yset\}\]
via alternating minimization or alternating gradient descent.
Two specific instances of the loss function $\loss(X,Y)$ arise from robust PCA and from the Gaussian factor model. 

\paragraph{Robust principal component analysis (RPCA)}

We study the robust PCA problem as formulated in \citet{agarwal2012noisy}, where the data matrix $Z\in\mathbb{S}^{n\times n}$ is generated from the model
\[ Z = \Aoper(\Xs + \Ys) + W,\]
where $\Aoper:\mathbb{S}^{d\times d} \rightarrow \mathbb{S}^{n\times n}$ is a linear operator mapping matrices from $\mathbb{S}^{d\times d}$ to $\mathbb{S}^{n\times n}$, and $W\in\mathbb{S}^{n\times n}$ represents a symmetric noise matrix. 

Our estimators are defined based on the least squares loss,
\begin{equation}\label{eqn:RPCA_program}(\Xh,\Yh)  = \argmin\left\{ \frac{1}{2}\fronorm{Z - \Aoper(X+Y)}^2: (X,Y)\in\Xset\times\Yset\right\}.\end{equation}
In determining the properties of the loss function, it is crucial that the operator $\Aoper$ satisfies certain conditions. Following the notion of RSC as introduced in \citet[Definition 2]{agarwal2012noisy}, we require $\Aoper$ to satisfy the following property:
\begin{assumption}\label{assump:RE} {\bf (Restricted Eigenvalue.)}
	There exist constants $\alpha_A,\beta_A$ and $\tau\geq 0$ such that for all $\Delta_X, \Delta_Y\in\R^{d\times d}$ with $\rank(\Delta_X)\leq 2r$,
	\[\alpha_A\left(\fronorm{\Delta_X}^2+\fronorm{\Delta_Y}^2\right) -\tau_{n,d} \leq\fronorm{\Aoper(\Delta_X+\Delta_Y)}^2 \leq  \beta_A\left(\fronorm{\Delta_X}^2+\fronorm{\Delta_Y}^2\right) +\tau_{n,d},\]
	where $\tau_{n,d}$ is given by
	\[\tau_{n,d} =  \tau\cdot \left(  \frac{\log d}{n^2} \norm{\Delta_Y}_1^2 + \sqrt{\frac{d^2\log d}{n^2}}\norm{\Delta_X}_\infty \norm{\Delta_Y}_1 \right).\]
\end{assumption}

Marginally in $X$ and in $Y$, by taking $\Delta_Y=0$ or $\Delta_X=0$, we see that  this assumption simply reduces to the well-known restricted eigenvalue property.
The expression $ \sqrt{\frac{d^2\log d}{n^2}}\norm{\Delta_X}_\infty \norm{\Delta_Y}_1$ reflects the restriction on the degree of interaction between $\Delta_X$ and $\Delta_Y$, which would hold if $\norm{\Aoper^{*}\Aoper(\Delta_X)}_\infty \approx \sqrt{\frac{d^2\log d}{n^2}}\norm{\Delta_X}_\infty$ ---for instance, an i.i.d. Gaussian ensemble will satisfy this property with high probability. We make this more general assumption on the operator $\Aoper$ to give a deterministic result on the least squares loss $\loss(X,Y)$. 

Now let the radii $\rho_X,\rho_Y$ satisfy $\rho_X, \rho_Y \leq c_0\cdot \sigma_r(\Xs)\kappa^{-1}(\Aoper)$ for some $c_0>0$, where $\sigma_r(\Xs)$ is the smallest singular value of $\Xs$, and where $\kappa(\Aoper)=\beta_A / \alpha_A$, which we can think of as a restricted condition number of the linear operator $\Aoper$ (i.e.~characterizing
the action of $\Aoper$ restricted to low-rank and sparse matrices). Given the initialization point $(X_0,Y_0)$, denote the corresponding neighborhoods by $\Xset_0=\Xset\cap\ball_2(X_0,\rho_X), \Yset_0=\Yset\cap\ball_2(Y_0,\rho_Y)$, and further assume that both the underlying matrices $(\Xs,\Ys)$ and the global optimal $(\Xh,\Yh)$ belong to these local neighborhoods $\Xset_0\times\Yset_0$. With this setup, we then have the following guarantee:

\begin{lemma}\label{lem:RPCA}
	Suppose that the sample size is large enough to satisfy
	\begin{equation}\label{eqn:RPCA_sample}\frac{32\tau\cdot sd \log d}{n^2}\leq\alpha_A. \end{equation}
	Then, under the previously stated conditions, if $\opnorm{\Aoper^*(W)}\leq c_1\cdot\alpha_A\sigma_r(\Xs)$, the steps $(X_t,Y_t)_{t=1}^{\infty}$ produced by the alternating minimization algorithm~\eqnref{altmin} satisfy
	\begin{equation*}\fronorm{(X_t,Y_t)-(\Xh,\Yh)} \leq   \overbrace{\left(1-\frac{\kappa^{-1}(\Aoper)}{3} \right)^{t}}^{\textnormal{linear convergence}}\cdot c_0\sigma_r(\Xs)\sqrt{\kappa^{-1}(\Aoper)} + \underbrace{c_2\cdot C\left(\fronorm{\Yh-\Ys}^2 + \frac{\spike^2}{\alpha_A^2}\frac{sd\log d}{n^2} \right)}_{\textnormal{statistical error term}}.\end{equation*}
	Here, $\{c_i >0,i=1,2\}$ are universal constants, and $C>0$ is defined in \thmref{altmin_contract}.
\end{lemma}
The proof of this lemma appears in \appref{RPCA_proofs}.

The result given in the lemma is the bound obtained by updating the $Y$ variable first instead of the $X$ variable. The statistical error involves the term $\spike^2\frac{sd\log d}{n^2}$, which appears as a consequence of the nonidentifiaibility of the model---see \citet{agarwal2012noisy} for a detailed discussion on the nonidentifiability of the matrix decomposition problem. With more effort, we can also prove the assumptions of \lemsref{altmin_inexact} and \lemssref{grad_descent} for the least square loss, implying that each step of the alternating minimization update can be replaced by several steps of the gradient descent updates. 

Several works also study alternating minimization and its variants such as alternating gradient descent employed on the robust PCA problem---for instance, see \citet{chen2015fast,yi2016fast,gu2016low}. In contrast to our approach, these methods are all based on the factorized approach, working on a highly nonconvex loss function arising from the factorization of the low-rank matrix $X=UU^\top$ (in the case of positive semidefinite matrix), which is identifiable only up to rotations. In fact, our framework can cover this factorization scenario as well, where the restricted strong convexity~\eqnref{RSC_obj} \& restricted smoothness~\eqnref{RSM_obj} of the original loss $\loss(X,Y)$ can be shown to be inherited by the reparametrized loss $\loss(UU^\top,Y)$, up to rotations, as long as the initialization condition \eqnref{initialization2} holds in the neighborhood $\Xset_0$. Although this approach is not the focus of our current paper, we refer the reader to \citet{chen2015fast}, \citet{tu2015low} and \citet{zheng2015convergent} for the related topic; see also \citet{gu2016low} and the references therein.

\paragraph{Gaussian factor model}

We next consider a Gaussian factor model, where our data consists of observations
\[ z_i = \Us w_i + \epsilon_i \;\text{ for $i=1,\ldots,n$}.\]
Here $\Us\in\R^{d\times r}$ represents the latent structure present in the data, while the other
terms in the model are the random factors $w_i\iidsim \mathcal{N}(0,\ident_r)$ and the independent noise $\epsilon_i\iidsim\mathcal{N}(0,\Ys)$. 

In the simplest setting, the covariance matrix of the noise, $\Ys$, is assume to be proportional to identity matrix in which case the dependence of the random vector $z_i$ is determined entirely by the latent structure $\Us$. Here we allow more flexibility on the covariance structure and assume that the covariance structure of the noise, $\Ys$, is sparse. Under this assumption, we can calculate $\Sigmas = \Us\Us{}^\top + \Ys$,  a low-rank + sparse decomposition, and can then estimate the unknown components $\Xs=\Us\Us{}^\top$ and $\Ys$ by solving the constrained optimization problem
\begin{equation}\label{eqn:factor_program}
(\Xh,\Yh)  =  \argmin \{\loss(X,Y):(X,Y)\in\Xset\times\Yset \} \text{ where }\loss(X,Y) = \inner{S_n}{(X+Y)^{-1}} - \log \det (X+Y)^{-1},\end{equation}
for $S_n=  \frac{1}{n}\sum_{i=1}^{n} z_i z_i^\top$, the sample covariance matrix of $z_i$'s,
where $\Xset$ and $\Yset$ are defined as in~\eqnref{XY_lowrank_sparse}. Note that the loss function~\eqnref{factor_program} is highly nonconvex due to the presence of the matrix inverse.

\citet{zwiernik2016maximum} study the related loss function $\loss(\Sigma) = \loss(X+Y)$ in the context of a linear Gaussian covariance model. The authors prove that this loss is in fact convex in the region\footnote{Specifically they show that the Hessian matrix of the loss function $\loss(\Sigma)$ is positive semidefinite in the region  $\{\Sigma\in\R^{d\times d}:0 \prec \Sigma \prec 2 S_n\}$.} $\{\Sigma\in\R^{d\times d}:0 \prec \Sigma \prec 2 S_n\}$ and further show that this region contains both the true covariance matrix $\Sigmas$ and the maximum likelihood estimator $\Sigmah$ with high probability, as long as the sample size is large enough, $n\gtrsim d$. In this regard, our setting can be seen as imposing different structure on the covariance matrix.

In the lemma to follow, we verify analogous results as in \citet{zwiernik2016maximum}, showing that the loss~\eqnref{factor_program} satisfies all the assumptions of \thmref{altmin_contract} in the local region, ensuring fast convergence of the alternating minimization algorithm. Suppose that the algorithm is initialized at the point $(X_0,Y_0)$ with the corresponding neighborhoods $\Xset_0=\Xset\cap\ball_2(X_0,\rho_X), \Yset_0=\Yset\cap\ball_2(Y_0,\rho_Y)$, where for some $c_0>0$ the radii are defined to satisfy
\[\rho_X,\rho_Y \leq c_0\cdot \min\{\sigma_r(\Xs) \kappa^{-3}(\Sigmas),\lammin(\Sigmas) \kappa^{-4}(\Sigmas) \}, \]
where $\sigma_r(\Xs)$ is the smallest singular value of $\Xs$, $\lammin(\Sigmas)$ (and $\lammax(\Sigmas)$ resp.) is the minimum (and maximum resp.) eigenvalue of $\Sigmas$, and where $\kappa(\Sigmas)=\lammax(\Sigmas)/\lammin(\Sigmas)$ is the condition number of $\Sigmas$. Assume also that these neighborhoods $\Xset_0\times\Yset_0$ contain the pair of true matrices $(\Xs,\Ys)$ and the global optimal $(\Xh,\Yh)$, i.e. $(\Xs,\Ys),(\Xh,\Yh)\in\Xset_0\times\Yset_0$. With this setup, we now establish the following probabilistic guarantee:
\begin{lemma}\label{lem:factor_model}
	Suppose that 
	\begin{equation}\label{eqn:factor_sample}\sqrt{\frac{d}{n}}\leq c_1\cdot \min\{\sigma_r(\Xs)\lammin^{-1}(\Sigmas)\kappa^{-4}(\Sigmas),\kappa^{-1}(\Sigmas) \}.\end{equation}
	Then, under the previously stated conditions, with probability at least $1-2e^{-d}$, the steps $(X_t,Y_t)_{t=1}^{\infty}$ produced by the alternating minimization~\eqnref{altmin} satisfy
	\begin{multline*}
	\fronorm{(X_t,Y_t)-(\Xh,\Yh)} \leq  \overbrace{\left(1-  c_2\kappa^{-4}(\Sigmas)\right)^{t}}^{\textnormal{linear convergence}} \cdot   c_3\min\{\sigma_r(\Xs)\kappa^{-1}(\Sigmas),\lammin(\Sigmas)\kappa^{-2}(\Sigmas)\} \\+ \underbrace{c_4\cdot C\left(\fronorm{\Yh-\Ys}^2 + \spike^2\frac{s}{d} \right)}_{\textnormal{statistical error term}}.
	\end{multline*}
	Here, $\{c_i >0,i=1,\ldots,4\}$ are universal constants, and $C>0$ is defined in \thmref{altmin_contract}. 
\end{lemma}
This lemma is proved in~\appref{factor_proofs}.

The discussion following~\lemref{RPCA} is also valid in this setting---in particular, the error due to the nonidentifiability of the model now appears as the term $\spike^2\frac{s}{d}$.

\subsection{Multitask regression}\label{sec:multitask}

Next we consider a class of regression problems where the response contains more than a single observation, that is, it takes multiple output values. Each component of the output corresponds to a specific task. 

Suppose we are given $m$ different tasks, where for each observation the response is of the form $z_i\in\R^m$. In the multitask regression model, the $m$ tasks are assumed to share the same feature vector, which we denote by $\phi_i\in\R^d$, and the response is generated through the linear model 
\begin{equation}\label{eqn:multitask_model}z_i = \Xs \phi_i + \epsilon_i,\end{equation}
where $\Xs\in\R^{m\times d}$ is an unknown matrix whose rows correspond to the underlying coefficient vectors for each task, and $\epsilon_i\in\R^m$ is the measurement error from a centered multivariate normal distribution, with an unknown covariance matrix $\textnormal{Cov}(\epsilon_i)=\Thetas{}^{-1}$. Given such model, we would like to recover the unknown matrices $\Xs$ and $\Thetas$ from the data $(z_i,\phi_i)_{i=1}^{n}$.
In the high-dimensional setting, it is common to estimate the matrix $\Xs$ under a low-rank constraint; this method is also referred to as ``reduced rank regression'' in the literature, e.g.  \citet{izenman1975reduced}. 

We would then like to optimize the constrained negative log-likelihood function,
\begin{equation}\label{eqn:multitask_program}
(\Xh, \Thetah)= \argmin\left\{ -\log\det(\Theta) + \frac{1}{n}\sum_{i=1}^{n} (z_i - X\phi_i)^\top \Theta (z_i- X\phi_i):X\in\Xset, \Theta\succeq 0\right\},
\end{equation}
where $\Xset=\{X\in\R^{m\times d}:\rank(X)\leq \rank(\Xs)=r\}$ represents the rank constraint on the coefficients $X$.\footnote{It is also possible to consider structural constraints on $\Theta$ or $\Sigma=\Theta^{-1}$ such as ``sparsity'' or ``low-rank + diagonal'' structure. For simplicity, we don't pursue this direction further, but our framework can be also applied to this general setting.}   A challenge here is that this problem is nonconvex in $(X,\Theta)$, in addition to the nonconvex constraint $X\in\Xset$, so it cannot be easily solved by standard convex optimization methods. Nevertheless, we will show that this problem satisfies all the assumptions that we require for our results, which allows to apply the (inexact) alternating minimization algorithm to solve the problem efficiently.

For the purpose of our analysis, we consider a random design model, i.e. the feature vectors are sampled from a Gaussian distribution $\phi_i\iidsim \mathcal{N}(0,\Sigma_\phi)$. Let $(X_0,\Theta_0)$ be the initialization point, and denote the local neighborhoods of the constraint sets around $(X_0,\Theta_0)$ by $\Xset_0 = \Xset\cap \ball_2(X_0,\rho_X)$ and $\Thetaset_0 = \mathbb{S}_{+}^{m\times m}\cap\ball_2(\Theta_0,\rho_\Theta)$. 
We choose the radii  $\rho_X, \rho_\Theta$ to satisfy $\rho_X\leq c_0\cdot  \sigma_r(\Xs)\kappa^{-1}(\Thetas)\kappa^{-1}(\Sigma_\phi)$ and $\rho_\Theta\leq c_0\cdot \lammin(\Thetas) \kappa^{-1}(\Sigma_\phi)$
for some $c_0>0$, where $\sigma_r(\Xs)$ is the smallest singular value of $\Xs$, $\lammin(\Thetas)$ is the smallest eigenvalue of $\Thetas$, and where $\kappa(\Thetas), \kappa(\Sigma_\phi)$ are the condition numbers of $\Thetas$ and $\Sigma_\phi$, respectively. Assume also that the initialization point $(X_0,\Theta_0)$ lies within these radii $\rho_X, \rho_\Theta$ to the unknown matrices $(\Xs,\Thetas)$ and the global optimal $(\Xh,\Thetah)$, i.e. $(\Xs,\Thetas), (\Xh,\Thetah)\in\Xset_0\times\Thetaset_0$. 

With these definitions in place, we have the following probabilistic guarantee:
\begin{lemma}\label{lem:multitask_regression}
	Suppose that 
	\begin{equation}\label{eqn:multitask_sample}
	\sqrt{\frac{1}{\lammin(\Thetas) \lammax(\Sigma_\phi)}} \sqrt{\frac{m+d}{n}}\leq c_1\cdot \sigma_r(\Xs)\kappa^{-1}(\Thetas)\kappa^{-1}(\Sigma_\phi).\end{equation}
	Then, under the previously stated conditions, with probability at least 
	$1- c_2\exp(-c_3(m+d))$, the steps $(X_t,\Theta_t)_{t=1}^{\infty}$ produced by the alternating minimization algorithm~\eqnref{altmin} satisfy 
	\begin{multline*}
	\fronorm{(X_t,\Theta_t) - (\Xh,\Thetah)} \leq \overbrace{\Big(1 - c_4(\kappa^{-1}(\Thetas)\kappa^{-1}(\Sigma_\phi) + \kappa^{-2}(\Thetas) ) \Big)^t }^{\textnormal{linear convergence}}\cdot (\textnormal{Const})
	\\+\underbrace{ c_5\cdot C\left(\fronorm{\Xh-\Xs}^2 +\frac{r(m+d)}{n}\frac{1}{\lammin(\Thetas) \lammax(\Sigma_\phi)}\right)}_{\textnormal{statistical error term}}
	\end{multline*}
	for all $t\geq 1$, where $(\textnormal{Const})$ is given by
	\[ (\textnormal{Const})= c_6\sigma_r(\Xs)\sqrt{\kappa^{-1}(\Thetas)\kappa^{-1}(\Sigma_\phi)}\cdot  \min\left\{1,\lammin^3(\Thetas)\kappa^2(\Thetas)\lammin(\Sigma_\phi)\right\}.\]
	Here, $\{c_i >0,i=1,\ldots,6\}$ are universal constants, and $C>0$ is  defined in \thmref{altmin_contract}.
\end{lemma}

By working with the inexact versions \lemsref{altmin_inexact} and \lemssref{grad_descent}, we can similarly obtain a linear rate of convergence for the alternating method when the minimization step for $X$ is approximated by successive iterates of gradient descent. 
(The alternating minimization update for $\Theta$ has a closed form solution, $\Theta=\left(\frac{1}{n}\sum_{i=1}^{n}(z_i - X\phi_i)(z_i - X\phi_i)^\top\right)^{-1}$, i.e.~the inverse of the sample covariance matrix.)
The result in \lemref{multitask_regression} assumes updating $\Theta$ first, but the analogous result is also available if the algorithm begins with $X$ updates instead.

Note that one can alternatively estimate $\Xs$ is by treating each task as a separate regression problem and simply performing least square procedure. However, it is generally known that the solution resulting from \eqnref{multitask_program} is significantly better than the least squares estimators, because the estimator defined in \eqnref{multitask_program} exploits the correlation structure  across tasks, whereas the least squares estimators ignore such correlation. We refer the reader to \citet{jain2015alternating} for a discussion of the similar type of results under the context of the pooled model.

\section{Empirical results}\label{sec:empirical}

We perform a numerical experiment on the multitask regression problem (\secref{multitask}) to examine the empirical performance of the alternating algorithm, as compared to performing gradient descent when treating $(x,y)$ as a single variable.\footnote{Code available at \url{http://www.stat.uchicago.edu/~rina/code/altmin_simulation.R}.} Fix the number of tasks $m=20$, the dimension of features $d=50$, and set the low-rank component $\Xs=U^{\star}V^{\star}{}^\top$ for rank $r=3$, where $U^{\star}\in\R^{20\times 3}$ and $V^{\star}\in\R^{50\times 3}$ are orthonormal matrices drawn uniformly at random. The features $\phi_i$ are drawn i.i.d. from the Gaussian distribution $\phi_i\iidsim\mathcal{N}(0,\Sigma_\phi)$, 
and the noise terms $\epsilon_i$ are generated as $\epsilon_i\iidsim\mathcal{N}(0,\Thetas{}^{-1})$, where $\Sigma_\phi$ and $\Thetas{}^{-1}$ are both defined to 
have a tapered covariance structure: we set $\Sigma_{\phi,ij} = 0.3^{|i-j|}$ and $\Thetas_{ij}{}^{-1} = \sigma^2\cdot \rho^{|i-j|}$, where $\rho$ is a local correlation parameter
that we vary, while $\sigma^2 =  \frac{\textnormal{Mean}(\fronorm{\Xs\phi_i}^2/m)}{3}$ is chosen to obtain a moderately difficult signal-to-noise ratio.
The responses, $z_i$, are then drawn according to the model \eqnref{multitask_model}.

The parameter $\rho$  controls the strength of the correlation of the noise (i.e.~correlation among entries of $\eps_i$, for a single observation $i$,
across the $m=20$ tasks). By varying $\rho$, we can vary the relative  condition numbers of the loss function $\loss(X,\Theta)$ given in~\eqnref{multitask_program} with
respect to the  variables $X$ and $\Theta$, i.e. $\kappa_X(\loss)$ versus $\kappa_{\Theta}(\loss)$. As discussed in \secref{conditioning},
convergence rates for alternating minimization type methods are expected to scale with the {\em minimum} of these two condition numbers,
while non-alternating methods (i.e.~gradient descent in the joint variable $(X,\Theta)$) will scale with the {\em maximum} of the two.

Given the data $(\phi_i,z_i)_{i=1}^{n}$ with sample size $n=200$, we solve the constrained minimization problem \eqnref{multitask_program} based on two 
iterative methods:
\begin{itemize}
	\item[$\bullet$] The alternating method, which alternates between updating $X$ and $\Theta$ at every iteration. For the $X$ update, fixing $\Theta$ we approximately
	minimize $\loss(X,\Theta)$ by taking one gradient descent step, while for the $\Theta$ update, fixing $X$ we minimize $\loss(X,\Theta)$ exactly:
	\[\begin{cases} &X_t = \Pr{\{\rank(X)\leq r\}}{X_{t-1} + \eta_X \cdot 2\Theta_{t-1} \left(\frac{1}{n}\sum_{i=1}^n (z_i - X_{t-1}\phi_i)\phi_i^\top\right)},\\ 
	& \Theta_t = \argmin_{\Theta\succeq 0}\loss(X_t,\Theta) = \left(\frac{1}{n}\sum_{i=1}^{n}(z_i - X_t\phi_i)(z_i - X_t\phi_i)^\top\right)^{-1}, \end{cases}\]
	with step size $\eta_X = 0.001$. \smallskip
	\item[$\bullet$] The joint gradient method, where we take gradient descent steps in the joint variable $(X,\Theta)$. The update step is given by
	\[\begin{cases}
	&X_t =  \Pr{\{\rank(X)\leq r\}}{X_{t-1} + \eta_X \cdot 2\Theta_{t-1} \left(\frac{1}{n}\sum_{i=1}^n (z_i - X_{t-1}\phi_i)\phi_i^\top\right) },\\
	&\Theta_t =  \Pr{\{\Theta\succeq 0\}}{\Theta_{t-1} + \eta_{\Theta}\cdot \left(\Theta_{t-1}^{-1} - \frac{1}{n}\sum_{i=1}^n (z_i - X_{t-1}\phi_i)(z_i - X_{t-1}\phi_i)^\top\right)},\end{cases}\]
	where we allow different step sizes on the two variables $X$ and $\Theta$.
	We set $\eta_X = 0.001$ as for the alternating method, 
	and select $\eta_\Theta\in\{\eta_1,\dots,\eta_{30}\}$, where $\eta_1,\dots,\eta_{30}$ is a geometric sequence from $\eta_1=5$ to $\eta_{30}=400$.
	For each trial, we then retain only the step size $\eta_\Theta$ that yields the lowest loss over any iteration, $\min_{t=1,\dots,T} \loss(X_t,\Theta_t)$,
	for the first $T=1200$ iterations. 
\end{itemize}

\begin{figure}[!t]\centering
	\includegraphics[width=0.33\textwidth]{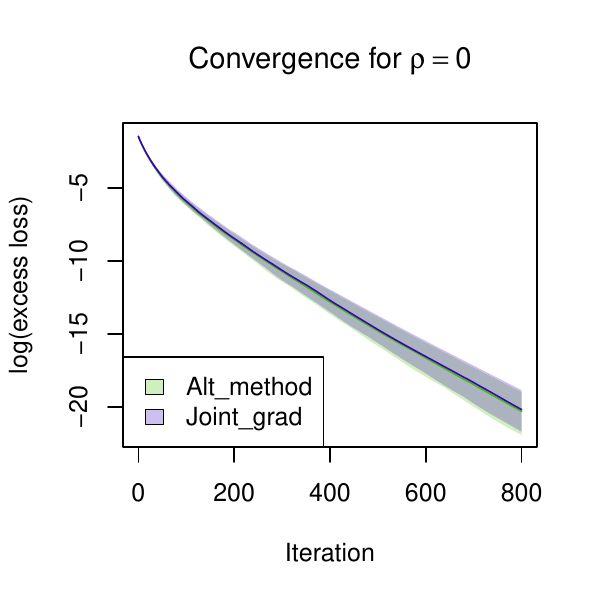}
	\includegraphics[width=0.33\textwidth]{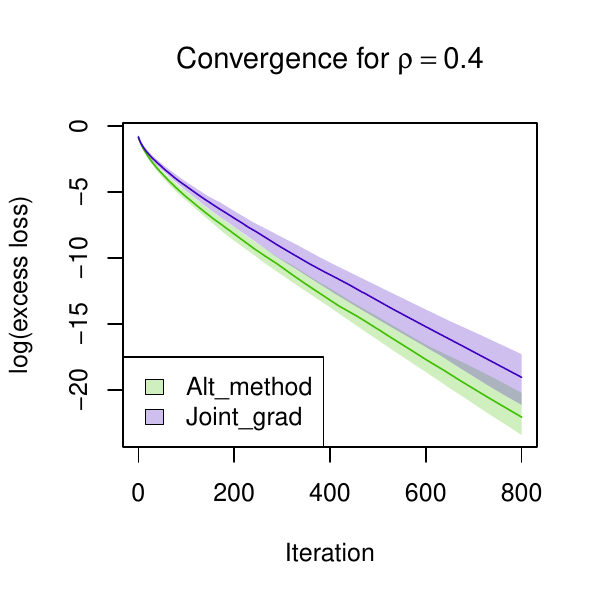}
	\includegraphics[width=0.33\textwidth]{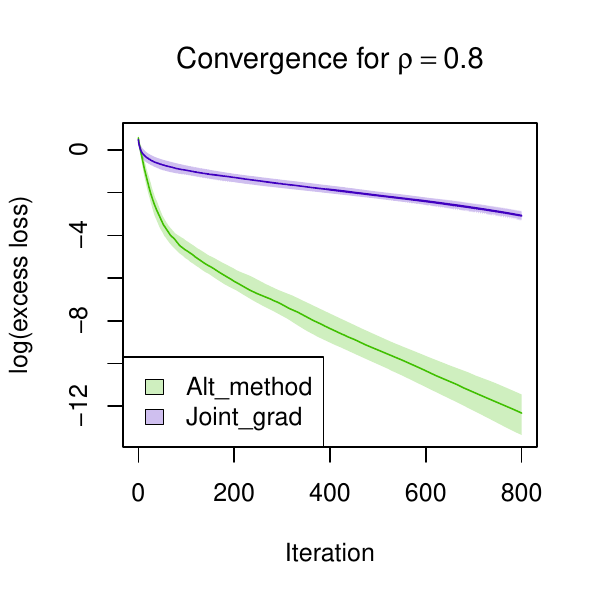}
	\caption{Comparison of the alternating algorithm and the joint gradient descent method applied to the simulated multitask regression problem. Both algorithms are initialized at $X_0=\Pr{\rank(X)\leq r}{X^{\textnormal{LS}}}$, where $X^{\textnormal{LS}}$ is the least squares estimator  ignoring the correlation structure.
	 Results are shown for three
		settings of the noise correlation parameter $\rho$, across iterations $t=0,1,\dots,600$. In each plot, the solid line indicates the median loss over $100$ trials, and 
		the light band shows the interquartile range. }
	\label{fig:multitask_simulation}
\end{figure}

For both method, we initialize at $X_0=\Pr{\rank(X)\leq r}{X^{\textnormal{LS}}}$, where $X^{\textnormal{LS}}$ is the least squares estimator  given by 
\[X^{\textnormal{LS},\top}= \left( \frac{1}{n}\sum_{i=1}^{n}\phi_i\phi_i^\top\right)^{-1}\cdot \left(\frac{1}{n}\sum_{i=1}^{n}\phi_i z_i^\top \right). \]

\figref{multitask_simulation} shows the excess loss at each iteration (on a log scale), where the excess loss is given by
\[\loss(X_t,\Theta_t)-\loss_{\textnormal{min}},\]
where $\loss_{\textnormal{min}}$ is the minimum loss achieved by either method over $T=1200$ iterations (calculated for each individual trial and each choice of $\rho$).
As clearly seen in the figure, for both methods, the errors scale linearly with the iteration number. Furthermore, comparing the two methods, we see that  they perform nearly identically 
when there is no correlation in the noise, i.e.~$\rho=0$; for low correlation, $\rho=0.4$, the alternating method is moderately faster,\footnote{Note that the 
	shaded bands in the plots are {\em not} standard error bars, but rather interquartile range over $100$ trials, so the difference between the two lines is indeed significant.} 
and for high correlation, $\rho=0.6$, the alternating method still shows rapid linear convergence while  joint gradient descent does not appear to converge well.
This is consistent with our theoretical results, since the alternating method scales with the better of the two condition numbers, i.e.~$\min\{\kappa_X(\loss),\kappa_{\Theta}(\loss)\}$,
while the joint gradient descent method is known to scale with the maximum. Since $\kappa_X(\loss)\sim \kappa(\Sigma_\phi)\kappa(\Thetas)$ 
while  $\kappa_{\Theta}(\loss)\sim \kappa^2(\Theta^\star)$, and $\kappa(\Sigma_\phi)$ is constant with respect to $\rho$ while $\kappa(\Theta)$
increases as the noise correlation parameter $\rho$ increases,
we see that the {\em minimum} condition number is less affected by increasing $\rho$, than the {\em maximum} condition number.

\section{Proofs of theorems}\label{sec:proofs_thms}
In this section, we prove our main result on linear convergence for the exact alternating minimization algorithm, \thmref{altmin_contract},
and for the inexact algorithm, \thmref{altmin_inexact}.
All other results are proved in the Appendix.

\subsection{Proof of \thmref{altmin_contract}}\label{sec:proof_thm_exact}

First we prove the bound on the $x$ update step, given in~\eqnref{thm1_x}, for iteration number $t$.
By definition of $x_t$, we can apply the first-order optimality condition~\eqnref{firstorder}, with $\Xset_0$ in place of $\Xset$,
to obtain
\[\inner{\xh - x_t}{\nabla_x \loss(x_t,y_{t-1})}\geq -\gamma(\Xset_0)\normx{\nabla_x \loss(x_t,y_{t-1})}^*\norm{x_t-\xh}^2_2.\]
Meanwhile, since $\xh$ is the minimizer of the problem $\min\{\loss(x,\yh):x\in\Xset_0\}$, we also have
\[ \inner{x_t-\xh }{\nabla_x\loss(\xh,\yh)} \geq -\gamma(\Xset_0)\normx{\nabla_x\loss(\xh,\yh)}^*\norm{x_t-\xh}^2_2.\]

Adding these two inequalities together, applying the initialization condition (\assumpref{initialization}) and rearranging terms several times, we have
\begin{align*}
&\frac{\alpha_x - \mu_x}{2} \norm{x_t -\xh}^2_2\\
&\geq\inner{x_t - \xh}{\nabla_x\loss(x_t,y_{t-1}) - \nabla_x\loss(\xh,\yh)} \\
&=\frac{1}{2}\biginner{\left(\begin{array}{c}x_t - \xh\\ y_{t-1} - \yh\end{array}\right)}{\nabla\loss(x_t,y_{t-1}) - \nabla\loss(\xh,\yh)} \\
&\hspace{0.4in}
+ \frac{1}{2}\inner{x_t - \xh}{\nabla_x\loss(x_t,y_{t-1}) - \nabla_x\loss(\xh,\yh)} 
- \frac{1}{2}\inner{y_{t-1}-\yh}{\nabla_y\loss(x_t,y_{t-1}) - \nabla_y\loss(\xh,\yh)}\\
&=\frac{1}{2}\biginner{\left(\begin{array}{c}x_t - \xh\\ y_{t-1} - \yh\end{array}\right)}{\nabla\loss(x_t,y_{t-1}) - \nabla\loss(\xh,\yh)} \\
&\hspace{0.4in}
+ \frac{1}{2}\inner{x_t - \xh}{\nabla_x\loss(x_t,\yh) - \nabla_x\loss(\xh,\yh)} 
- \frac{1}{2}\inner{y_{t-1}-\yh}{\nabla_y\loss(\xh,y_{t-1}) - \nabla_y\loss(\xh,\yh)}\\
&\hspace{0.8in}
+\frac{1}{2}\Big[\inner{x_t - \xh}{\nabla_x \loss(x_t,y_{t-1}) -\nabla_x\loss(x_t,\yh)}
- \inner{y_{t-1}-\yh}{\nabla_y\loss(x_t,y_{t-1}) - \nabla_y\loss(\xh,y_{t-1})}\Big]\\
&\geq\frac{\alpha_x}{2} \norm{x_t-\xh}_2^2 + \frac{\alpha_y}{2}\norm{y_{t-1}-\yh}_2^2 - \frac{\alpha_x\errx^2+\alpha_y\erry^2 }{2}\\
&\hspace{0.4in}
+ \frac{1}{2}\inner{x_t - \xh}{\nabla_x\loss(x_t,\yh) - \nabla_x\loss(\xh,\yh)} 
- \frac{1}{2}\inner{y_{t-1}-\yh}{\nabla_y\loss(\xh,y_{t-1}) - \nabla_y\loss(\xh,\yh)}\\
&\hspace{0.8in} - \frac{1}{2}\left(\frac{1}{2}\mu_x \norm{x_t-\xh}_2^2 + \frac{1}{2}\mu_y \norm{y_{t-1}-\yh}_2^2 + \alpha_x\errx^2+\alpha_y\erry^2\right),
\end{align*}
where the last step holds by applying joint restricted strong convexity (\assumpref{RSC_joint})
to the first term, and the cross-product condition (\assumpref{crossprod}) to the 
expression in square brackets. (Note that these assumptions can be applied since we have $x_t\in\Xset_0$ and $y_{t-1}\in\Yset_0$). Combining terms and simplifying,  multiplying by $2$, and using the assumption that $\mu_x\geq 0$ while $\mu_y\leq \alpha_y$,
we obtain
\begin{equation}\label{eqn:contract_rearrange}
0 \geq \diver{x_t}{\xh} - \diver{y_{t-1}}{\yh} + \frac{\alpha_y}{2}\norm{y_{t-1}-\yh}^2_2 - 2\alpha_x\errx^2  -  2\alpha_y\erry^2.
\end{equation}
Now, by restricted smoothness (\assumpref{RSM}) and using the assumption $\alpha_y\leq \beta_y$,
\[\frac{\alpha_y}{2}\norm{y_{t-1}-\yh}^2_2 \geq \frac{\alpha_y}{2\beta_y}\diver{y_{t-1}}{\yh}- \frac{\alpha_y}{2}\erry^2.\]
Returning to~\eqnref{contract_rearrange} and rearranging terms,
\[\diver{x_t}{\xh} \leq \left(1 - \frac{\alpha_y}{2\beta_y}\right)\diver{y_{t-1}}{\yh} + 3\left(\alpha_x\errx^2 + \alpha_y\erry^2 \right).\]
If $\diver{x_t}{\xh}\geq 0$, then by taking a square root on both sides, we obtain
\[\sdiver{x_t}{\xh}  \leq \sqrt{1-\frac{\alpha_y}{2\beta_y}}\cdot\sdiver{y_{t-1}}{\yh} + \sqrt{3(\alpha_x\errx^2+\alpha_y\erry^2)},\]
thus proving that the bound~\eqnref{thm1_x} holds at time $t$, while if $\diver{x_t}{\xh}\leq 0$, then $\sdiver{x_t}{\xh}=0$ and so the bound holds trivially.
The proof that the analogous bound~\eqnref{thm1_y} for the $y$ update step, proceeds similarly.

By applying these bounds recursively, along with the restricted strong convexity and restricted smoothness conditions,
we can obtain the result~\eqnref{converge_L2} showing linear convergence in the $\ell_2$ norm;
details are given in \appref{l2_details}.

\subsection{Proof of \thmref{altmin_inexact}}\label{sec:proof_altmin_inexact}

This proof is a straightforward combination of the relaxed triangle inequality (\assumpref{triangle}) with \thmref{altmin_contract}, the contraction
result for the exact alternating minimization algorithm.
First, since $x_t^{\textnormal{exact}}$ exactly solves the alternating minimization step, i.e.~$\argmin_{x\in\Xset_0}\loss(x,y_{t-1})$,
\thmref{altmin_contract} proves that
\[\sdiver{x_t^{\textnormal{exact}}}{\xh}\leq \sqrt{1 - \frac{\alpha_y}{2\beta_y}}\sdiver{y_{t-1}}{\yh} + \sqrt{3(\alpha_x\errx^2+\alpha_y\erry^2)}.\]
Next, we use this to bound $\sdiver{x_t}{\xh}$, using only the assumption that $x_t$ is chosen to be within radius $\delta_t^x$ of $x_t^{\textnormal{exact}}$.
By the relaxed triangle inequality (\assumpref{triangle}), 
\begin{align*}
\sdiver{x_t}{\xh} 
&\leq \sdiver{x_t^{\textnormal{exact}}}{\xh} + \sqrt{\beta_x}\norm{x_t - x_t^{\textnormal{exact}}}_2 + \sqrt{\alpha_x}\errx\\
&\leq \left(\sqrt{1 - \frac{\alpha_y}{2\beta_y}}\cdot \sdiver{y_{t-1}}{\yh} + \sqrt{3(\alpha_x\errx^2+\alpha_y\erry^2)}\right) + \sqrt{\beta_x}\delta_t^x + \sqrt{\alpha_x}\errx\\
&\leq \sqrt{1 - \frac{\alpha_y}{2\beta_y}}\cdot \sdiver{y_{t-1}}{\yh} +  \sqrt{\beta_x}\delta_t^x+\sqrt{8(\alpha_x\errx^2+\alpha_y\erry^2)}.
\end{align*}
This proves the bound~\eqnref{thm2_x}.
The bound~\eqnref{thm2_y} on $\sdiver{y_t}{\yh}$ is proved analogously.

\section{Discussion}\label{sec:discussion}

In this paper, we present a general convergence result for the alternating minimization method and its inexact variants in the presence of nonconvex constraints. We have shown that under standard assumptions on the loss function $\loss(x,y)$, these methods can offer linear convergence, as long as certain initialization condition is satisfied. A major tool allowing to handle the nonconvex constraints is the local concavity coefficient, which enables us to bound the concavity arising from the nonconvexity of the sets. The important implication of our result is the computational gain inherent in the alternating methods that relies on the variable with better conditioned, in contrast to the methods that work on the combined variable. This phenomenon has been further demonstrated through numerical experiments.

There are many open questions related to our problem. First, the initialization condition presented in this paper requires the careful choice of initialization point to begin with, but in some applications, a valid initialization procedure may not be available. Characterizing the convergence behavior of alternating minimization 
without strong initialization assumptions is therefore of practical importance for real applications. Furthermore, as seen in the low-rank factorization approach, it would be interesting to generalize our current result to more general setup, including  loss function that may only satisfy convexity up to 
identifiability issues (such as matrix factorization methods), which would extend our understanding of alternating minimization and descent
algorithms to a  broader range of problems.

\subsection*{Acknowledgements} 
R.F.B.~was supported by an Alfred P.~Sloan Fellowship and by  NSF award DMS-1654076.

\bibliographystyle{plainnat}
\bibliography{altmin.bib}

\appendix


\section{Additional proofs}\label{app:add_proofs}

\subsection{Proof of $\ell_2$ convergence bound~\eqnref{converge_L2}}\label{app:l2_details}
Here we give details for the $\ell_2$ convergence bound~\eqnref{converge_L2} in \thmref{altmin_contract}.
We first write $r_x = \sqrt{1 - \frac{\alpha_y}{2\beta_y}}$ and $r_y = \sqrt{1 - \frac{\alpha_x}{2\beta_x}}$
for simplicity; then~\eqnref{thm1_x} and~\eqnref{thm1_y} can be rewritten as
\[\sdiver{x_t}{\xh}\leq r_x\sdiver{y_{t-1}}{\yh}+\sqrt{3(\alpha_x\errx^2+\alpha_y\erry^2)},\]
and
\[\sdiver{y_t}{\yh}\leq r_y\sdiver{x_t}{\xh} + \sqrt{3(\alpha_x\errx^2+\alpha_y\erry^2)}.\]
Then, 
applying these bounds recursively, we have
\[\sdiver{x_t}{\xh}\leq r_x(r_xr_y)^{t-1}\sdiver{y_0}{\yh} + \sqrt{3(\alpha_x\errx^2+\alpha_y\erry^2)}\cdot \frac{1+r_x}{1-r_xr_y},\]
and
\[\sdiver{y_t}{\yh}\leq (r_xr_y)^t\sdiver{y_0}{\yh} + \sqrt{3(\alpha_x\errx^2+\alpha_y\erry^2)}\cdot \frac{1+r_y}{1-r_xr_y}.\]
Let $\alpha_{\textnormal{min}}=\min\{\alpha_x,\alpha_y \}, \alpha_{\textnormal{max}}=\max\{\alpha_x,\alpha_y \}$. By joint restricted strong convexity~\eqnref{RSC_joint},
\begin{multline*}\norm{(x_t,y_t)-(\xh,\yh)}_2 = \sqrt{\norm{x_t-\xh}^2_2 + \norm{y_t - \yh}^2_2}\\ \leq \sqrt{\frac{(\diver{x_t}{\xh} + \diver{y_t}{\yh})}{\alpha_{\textnormal{min}}} + \frac{2\alpha_{\textnormal{max}}(\errx^2 + \erry^2)}{\alpha_{\textnormal{min}}}}
\leq \frac{\sdiver{x_t}{\xh} + \sdiver{y_t}{\yh}}{\sqrt{\alpha_{\textnormal{min}}}} + \frac{\sqrt{2\alpha_{\textnormal{max}}(\errx^2 + \erry^2)}}{\sqrt{\alpha_{\textnormal{min}}}}\\
\leq (r_xr_y)^t \cdot \sdiver{y_0}{\yh}\cdot \frac{(1+r_y^{-1})}{\sqrt{\alpha_{\textnormal{min}}}} + \frac{\sqrt{3(\alpha_x\errx^2+\alpha_y\erry^2)}\cdot \left(\frac{2+2r_x}{1-r_xr_y}\right)}{\sqrt{\alpha_{\textnormal{min}}}} +\frac{\sqrt{2\alpha_{\textnormal{min}}(\errx^2 + \erry^2)}}{\sqrt{\alpha_{\textnormal{min}}}}
\end{multline*}
and, by definition of $r_x$ and $r_y$ and the fact that $r_xr_y\leq 1$, we see that $r_x,r_y\in[1/\sqrt{2},\sqrt{2}]$. Simplifying,
\[\norm{(x_t,y_t)-(\xh,\yh)}_2
\leq (r_xr_y)^t \cdot \sdiver{y_0}{\yh}\cdot \frac{\sqrt{6}}{\sqrt{\alpha_{\textnormal{min}}}} + \frac{\sqrt{3(\alpha_x\errx^2+\alpha_y\erry^2)}\cdot\frac{2+2\sqrt{2}}{1-r_xr_y}}{\sqrt{\alpha_{\textnormal{min}}}} + \frac{\sqrt{2\alpha_{\textnormal{max}}(\errx^2 + \erry^2)}}{\sqrt{\alpha_{\textnormal{min}}}},\]
and by restricted smoothness~\eqnref{RSM_marg},
\[\diver{y_0}{\yh}\leq \beta_y\norm{y_0 - \yh}^2_2 + \alpha_y\erry^2.\]
Combining everything, and simplifying, we obtain the overall convergence guarantee~\eqnref{converge_L2}.

\subsection{Proof of \lemref{altmin_inexact}}\label{app:lem_altmin_inexact}

For convenience define
\[r_x = \sqrt{1 - \frac{\alpha_y}{2\beta_y}} + (1+\sqrt{2}) \cdot c_x\sqrt{\frac{\beta_x}{\alpha_x}}
\text{ and }r_y  = \sqrt{1 - \frac{\alpha_x}{2\beta_x}} + (1+\sqrt{2}) \cdot c_y\sqrt{\frac{\beta_y}{\alpha_y}}.\]
(Comparing to the proof of the $\ell_2$ convergence bound~\eqnref{converge_L2} for the exact algorithm, given in \appref{l2_details}, we see that these definitions coincide with the previous ones in the special
case that $c_x = c_y=0$, i.e.~when our updates are exact.) Define also $D_0 = \sqrt{\alpha_x}\rho_x + \sqrt{\beta_y}\rho_y$.

We will first show, by induction, that for each $t\geq 1$,
\begin{equation}\label{eqn:inexact_induct}
\begin{cases}
\sdiver{x_t}{\xh} \leq r_x \cdot (r_x r_y)^{t-1} \cdot D_0 + \frac{1+r_x}{1-r_xr_y}\cdot C' \max\{\errx,\erry\},\\
\sdiver{y_t}{\yh} \leq (r_x r_y)^t \cdot D_0 + \frac{1+r_y}{1-r_xr_y}\cdot C' \max\{\errx,\erry\},
\end{cases}
\end{equation}
where
\begin{equation}\label{eqn:Cprime}
C' =4\left(1+c_x\sqrt{\frac{\beta_x}{\alpha_x}} + c_y\sqrt{\frac{\beta_y}{\alpha_y}}\right)\sqrt{\alpha_x+\alpha_y} + C_x\sqrt{\beta_x} + C_y\sqrt{\beta_y}.\end{equation}

First we prove the bounds~\eqnref{inexact_induct} at time $t=1$.
For the $x$ bound,
\begin{align*}
&\sdiver{x_1}{\xh}
\leq \sqrt{1 - \frac{\alpha_y}{2\beta_y}}\cdot \sdiver{y_0}{\yh} +  \sqrt{\beta_x}\delta_t^x+ \sqrt{8(\alpha_x\errx^2+\alpha_y\erry^2)} \text{ by \thmref{altmin_inexact}}\\
&\leq \sqrt{1 - \frac{\alpha_y}{2\beta_y}}\cdot \sdiver{y_0}{\yh} +  \sqrt{\beta_x} \left( c_x\norm{x_0-x_1^{\textnormal{exact}}}_2 + C_x\errx\right)+ \sqrt{8(\alpha_x\errx^2+\alpha_y\erry^2)} \text{ by~\eqnref{delta_rule_recursive}}\\
&\leq \sqrt{1 - \frac{\alpha_y}{2\beta_y}}\cdot (\sqrt{\beta_y}\rho_y+\sqrt{\alpha_y}\erry) +  \sqrt{\beta_x} \left( c_x\cdot \rho_x+ C_x\errx\right)+ \sqrt{8(\alpha_x\errx^2+\alpha_y\erry^2)},
\end{align*}
where the last step holds since $x_1^{\textnormal{exact}}\in\Xset_0\subset\ball_2(x_0,\rho_x)$, and $\sdiver{y_0}{\yh}$ can be bounded by restricted smoothness (\assumpref{RSM}).
Simplifying,
\[\sdiver{x_1}{\xh} \leq r_x D_0 + C'\max\{\errx,\erry\},\]
which proves the bound~\eqnref{inexact_induct} on $\sdiver{x_1}{\xh}$ at time $t=1$.
Similarly, for the $y$ bound,
\begin{align*}
&\sdiver{y_1}{\yh}
\leq \sqrt{1 - \frac{\alpha_x}{2\beta_x}}\cdot \sdiver{x_1}{\xh} +  \sqrt{\beta_y}\delta_t^y+ \sqrt{8(\alpha_x\errx^2+\alpha_y\erry^2)} \text{ by \thmref{altmin_inexact}}\\
&\leq \sqrt{1 - \frac{\alpha_x}{2\beta_x}}\cdot \sdiver{x_1}{\xh} +  \sqrt{\beta_y} \left( c_y\norm{y_0-y_1^{\textnormal{exact}}}_2 + C_y\erry\right)+ \sqrt{8(\alpha_x\errx^2+\alpha_y\erry^2)} \text{ by~\eqnref{delta_rule_recursive}}\\
&\leq \sqrt{1 - \frac{\alpha_x}{2\beta_x}}\cdot \left(r_x D_0 + C'\max\{\errx,\erry\}\right) +  \sqrt{\beta_y} \left( c_y\rho_y + C_y\erry\right)+ \sqrt{8(\alpha_x\errx^2+\alpha_y\erry^2)}\\
&\leq r_x r_y D_0 + (1+r_y)\cdot C'\max\{\errx,\erry\},
\end{align*}
where for the last step we use the fact that $r_y\geq \sqrt{1-\frac{\alpha_x}{2\beta_x}}$ by definition.

Next, take any $t\geq 2$. For the $x$ bound, we first calculate
\begin{align*}
\norm{x_{t-1}-x_t^{\textnormal{exact}}}_2 
&\leq \norm{x_{t-1}-\xh}_2 + \norm{x_t^{\textnormal{exact}}-\xh}_2 \\
&\leq \frac{1}{\sqrt{\alpha_x}}\left(\sdiver{x_{t-1}}{\xh} + \sdiver{x_t^{\textnormal{exact}}}{\xh}\right) + \frac{2\sqrt{\alpha_{x}\errx^2 +\alpha_y \erry^2}}{\sqrt{\alpha_{x}}}\text{ by joint restricted strong convexity~\eqnref{RSC_joint}}\\
&\leq  \frac{1}{\sqrt{\alpha_{x}}}\left(\sdiver{x_{t-1}}{\xh} + \sdiver{y_{t-1}}{\yh} + \sqrt{3(\alpha_x\errx^2+\alpha_y\erry^2)}\right) +\frac{2\sqrt{\alpha_{x}\errx^2 + \alpha_y\erry^2}}{\sqrt{\alpha_{x}}}\text{ by \thmref{altmin_contract}}\\
&\leq  \frac{1}{\sqrt{\alpha_{x}}}\left(\sdiver{x_{t-1}}{\xh} + \sdiver{y_{t-1}}{\yh} + 4\sqrt{\alpha_{x}\errx^2+\alpha_y\erry^2}\right).
\end{align*}
We now bound $\sdiver{x_t}{\xh}$:
\begin{align*}
&\sdiver{x_t}{\xh}
\leq \sqrt{1 - \frac{\alpha_y}{2\beta_y}}\cdot \sdiver{y_{t-1}}{\yh} +  \sqrt{\beta_x}\delta_t^x+ \sqrt{8(\alpha_x\errx^2+\alpha_y\erry^2)} \text{ by \thmref{altmin_inexact}}\\
&\leq \sqrt{1 - \frac{\alpha_y}{2\beta_y}}\cdot \sdiver{y_{t-1}}{\yh} +  \sqrt{\beta_x} \left( c_x\norm{x_{t-1}-x_t^{\textnormal{exact}}}_2 + C_x\errx\right)+ \sqrt{8(\alpha_x\errx^2+\alpha_y\erry^2)} \text{ by~\eqnref{delta_rule_recursive}}\\
&\leq \sqrt{1 - \frac{\alpha_y}{2\beta_y}}\cdot \sdiver{y_{t-1}}{\yh} +  \sqrt{\beta_x} \left( c_x\left[\frac{1}{\sqrt{\alpha_{x}}}\left(\sdiver{x_{t-1}}{\xh} + \sdiver{y_{t-1}}{\yh} + 4\sqrt{\alpha_{x}\errx^2+\alpha_y\erry^2}\right) \right]+ C_x\errx\right)\\
&\hspace{1in}+ \sqrt{8(\alpha_x\errx^2+\alpha_y\erry^2)}\\
&\leq  \left(\sqrt{1 - \frac{\alpha_y}{2\beta_y}} + c_x\sqrt{\frac{\beta_x}{\alpha_x}}\right)\sdiver{y_{t-1}}{\yh} + c_x\sqrt{\frac{\beta_x}{\alpha_x}}\sdiver{x_{t-1}}{\xh} +C'\max\{\errx,\erry\},
\end{align*}
where $C'$ is defined as in~\eqnref{Cprime} above.
Assuming by induction that the bounds~\eqnref{inexact_induct} hold with $t-1$ in place of $t$, we obtain
\begin{multline*}
\sdiver{x_t}{\xh}\leq \left(\sqrt{1 - \frac{\alpha_y}{2\beta_y}} + c_x\sqrt{\frac{\beta_x}{\alpha_x}}\right)\cdot \left((r_x r_y)^{t-1} D_0 + \frac{1+r_y}{1-r_xr_y}\cdot C' \max\{\errx,\erry\}\right)\\ + c_x\sqrt{\frac{\beta_x}{\alpha_x}}\left(r_x(r_x r_y)^{t-2} D_0 + \frac{1+r_x}{1-r_xr_y}\cdot C' \max\{\errx,\erry\}\right)  +C'\max\{\errx,\erry\}.\end{multline*}
Since $r_y\geq 1/\sqrt{2}$ we can rewrite this as
\begin{multline*}
\sdiver{x_t}{\xh}\leq \left(\sqrt{1 - \frac{\alpha_y}{2\beta_y}} + c_x\sqrt{\frac{\beta_x}{\alpha_x}}\right)\cdot \left((r_x r_y)^{t-1} D_0 +  \frac{1+r_y}{1-r_xr_y}\cdot C'  \max\{\errx,\erry\}\right)\\ + c_x\sqrt{\frac{\beta_x}{\alpha_x}}\left(\sqrt{2}(r_x r_y)^{t-1} D_0 +  \frac{1+r_x}{1-r_xr_y}\cdot C'  \max\{\errx,\erry\}\right)  +C'\max\{\errx,\erry\}.\end{multline*}
Plugging in the definition of $r_x$, then,
\[\sdiver{x_t}{\xh}\leq r_x\cdot (r_xr_y)^{t-1} \cdot D_0 + \left[c_x\sqrt{\frac{\beta_x}{\alpha_x}} \cdot \frac{1+r_x}{1-r_xr_y} +\left(\sqrt{1 - \frac{\alpha_y}{2\beta_y}} + c_x\sqrt{\frac{\beta_x}{\alpha_x}}\right) \cdot\frac{1+r_y}{1-r_xr_y} + 1\right] \cdot C' \max\{\errx,\erry\}.\]
Plugging in the definition of $r_x$, and the assumption that $r_xr_y<1$, we see that the term in square brackets is bounded by $\frac{1+r_x}{1-r_xr_y}$, which proves 
the desired bound on $\sdiver{x_t}{\xh}$ as in~\eqnref{inexact_induct}, as desired. The bound on $\sdiver{y_t}{\yh}$ is proved similarly.

Finally, by joint restricted strong convexity~\eqnref{RSC_joint}, we know that
\[\norm{x_t - \xh}_2 \leq \frac{\sdiver{x_t}{\xh}}{\sqrt{\alpha_x}} + \frac{\sqrt{\alpha_x\errx^2+\alpha_y\erry^2}}{\sqrt{\alpha_x}} \text{ and }\norm{y_t - \yh}_2 \leq \frac{\sdiver{y_t}{\yh}}{\sqrt{\alpha_y}} +  \frac{\sqrt{\alpha_x\errx^2+\alpha_y\erry^2}}{\sqrt{\alpha_y}} .\]
Combining this with the bounds~\eqnref{inexact_induct} proves the result.

\subsection{Proof of \lemref{LCC_rho}}\label{app:lem_LCC_rho}

Take any $x,x' \in\Xset_0\subset \Xset$ and take $t\in[0,1]$.
By the curvature condition (\defref{curvature}) on the larger set $\Xset$, we can find a family
of points $\xt_t\in\Xset$, indexed by $t\in[0,1]$, such that $\delta_t\rightarrow 0$, where
\[ \normx{\big((1-t)x + tx'\big) - \xt_t} \leq t\cdot\Big[ \gamma_x(\Xset) \cdot  \norm{x-x'}^2_2 + \delta_t \Big]. \]
Next, we show that $\xt_t\in\Xset_0$ for sufficiently small $t>0$. Recall that $\Xset_0 = \Xset\cap \ball_2(x_0,\rho_x)$,
and therefore we only need to check that $\norm{\xt_t - x_0}_2\leq \rho_x$. 
Since $\norm{\cdot}_2\leq \normx{\cdot}$ by assumption,
we have
\begin{align*}
\norm{\xt_t - x_0}_2
&\leq \norm{\xt_t - \big((1-t)x+tx'\big)}_2 + \norm{\big((1-t)x+tx'\big) - x_0}_2 \\
&\leq \normx{\xt_t - \big((1-t)x+tx'\big)} + \norm{\big((1-t)x+tx'\big) - x_0}_2 \\
&\leq  t\cdot\Big[ \gamma_x(\Xset) \cdot  \norm{x-x'}^2_2 + \delta_t \Big]+ \norm{\big((1-t)x+tx'\big) - x_0}_2 .
\end{align*}
Next, a simple calculation shows that
\begin{align*}
\norm{\big((1-t)x+tx'\big) - x_0}_2
& = \norm{(1-t)\cdot(x-x_0) + t\cdot(x'-x_0)}_2\\
& =  \sqrt{(1-t)\norm{x-x_0}^2_2 + t\cdot \norm{x'-x_0}^2_2 -  t(1-t)\norm{x-x'}^2_2},
\end{align*}
and since $x,x'\in\Xset_0\subset\ball_2(x_0,\rho_x)$, we obtain
\[\norm{\big((1-t)x+tx'\big) - x_0}_2\leq \sqrt{\rho^2_x -  t(1-t)\norm{x-x'}^2_2}\leq\rho_x - \frac{t(1-t)\norm{x-x'}^2_2}{2\rho_x}.\]
Combining everything,
\[\norm{\xt_t - x_0}_2 \leq \rho_x - t\norm{x-x'}^2_2 \cdot \left[ \frac{1}{2\rho_x} - \gamma_x(\Xset) - \frac{t}{2\rho_x}- \frac{\delta_t}{\norm{x-x'}^2_2}\right].\]
Since $\gamma_x(\Xset)<\frac{1}{2\rho_x}$ by assumption, and $\delta_t\rightarrow0$, we can find some
$t_0>0$ such that, for all $t\in[0,t_0]$,
\[\frac{t}{2\rho_x}+\frac{\delta_t}{\norm{x-x'}^2_2} \leq \frac{1}{2\rho_x} - \gamma_x(\Xset).\]
Therefore, $\xt_t\in\Xset_0$ for all $t\in[0,t_0]$, and so
\[\frac{\min_{x''\in\Xset_0}\normx{\big((1-t)x + tx'\big) - x''} }{t} \leq \frac{\normx{\big((1-t)x + tx'\big)- \xt_t} }{t} \leq  \gamma_x(\Xset) \cdot  \norm{x-x'}^2_2 + \delta_t\]
for all $t\in[0,t_0]$. This proves that
\[\lim_{t\rightarrow 0}\frac{\min_{x''\in\Xset_0}\norm{\big((1-t)x + tx'\big) - x''} }{t} \leq \gamma_x(\Xset)\cdot  \norm{x-x'}^2_2.\]
Since $x,x'\in\Xset_0$ were chosen arbitrarily, then, we have shown that
\[\gamma_x(\Xset_0)\leq \gamma_x(\Xset).\]

%
%
%

\section{Proofs for examples (\secref{examples})}\label{app:examples_proofs}

In this section, we provide the proofs of the lemmas displayed in \secref{examples}. Throughout the section, given any function $\mathsf{f}(A)$ over a matrix variable $A\in\R^{m\times n}$, we   write $\nabla^2_{AA}\mathsf{f}(A)\in\R^{mn\times mn}$ to refer to the second derivative of $\mathsf{f}(A)$ with respect to the {\em vectorized} variable $\vect{A}\in\R^{mn}$.

\subsection{Proof of \lemref{matrix_lcc}}

We first reparametrize the variable $X\in\Xset$ by $X=\ufun(U)=UU^\top$ with the corresponding convex set
\[\Uset = \left\{U\in\R^{d\times r}: \max_{i=1,\ldots,d}\norm{U_{i*}}_2 \leq  \sqrt{\frac{\spike}{d}}  \right\},\]
where $U_{i*}$ represents $i$th row of $U$. Note that under such reparametrization, we trivially have $\Xset=\ufun(\Uset)$.
Now take $X,X'\in\Xset$ with $X=UU^\top, X'=U'U'^\top$. For $t>0$, let $X_t = (1-t)X + tX'$ and $U_t = (1-t)U + tU'$. Then, by Taylor's theorem,
\begin{align}\label{eqn:reparam_curv}
&X_t  - \ufun(U_t) = (1-t) \ufun(U) + t \ufun(U') - \ufun(U_t)  \nonumber\\
&= (1-t) (\ufun(U) - \ufun(U_t)) + t (\ufun(U') - \ufun(U_t))\nonumber \\
&= (1-t)\left[\nabla \ufun(U_t) (U - U_t) + \frac{1}{2}\nabla^2 \ufun(U_*) (U - U_t, U - U_t) \right] + t\left[\nabla \ufun(U_t) (U' - U_t) + \frac{1}{2}\nabla^2 \ufun(U_\#) (U' - U_t, U' - U_t) \right] \nonumber\\
&= (1-t)\left[t \nabla \ufun(U_t) (U-U') + \frac{t^2}{2}\nabla^2 \ufun(U_*) (U-U', U-U') \right]\nonumber\\
&\hspace{2in} + t\left[(1-t)\nabla \ufun(U_t) (U' - U) + \frac{(1-t)^2}{2}\nabla^2 \ufun(U_\#) (U' - U, U'- U) \right].
\end{align}
Meanwhile, some calculation yields that for $i, j = 1,\ldots,d$, 
\[\nabla^2 \ufun_{ij}(U) = (e_i e_j^\top \otimes \ident_r + e_j e_i^\top \otimes \ident_r) \in\R^{dr \times dr},\]
where $e_i\in\R^{d}$ denotes the $i$th standard basis vector. Hence, we have
\[\nabla^2 \ufun(U_*) (U-U', U-U')=\nabla^2 \ufun(U_\#) (U-U', U-U') = 2(U-U')(U-U')^\top.  \]
Combining with \eqnref{reparam_curv}, 
\[\min_{X''\in\Xset}\nucnorm{X'' - X_t} \leq \nucnorm{\ufun(U_t) - X_t} = t(1-t)\nucnorm{(U-U')(U-U')^\top}= t(1-t)\fronorm{U-U'}^2,\]
so dividing out by $t$ and taking $t\to 0$,
\[\limsup_{t\to 0}\frac{\min_{X''\in\Xset}\normx{X'' - X_t} }{t} \leq\fronorm{U-U'}^2 \leq \frac{5}{4\sigma_r(X)}\fronorm{X-X'}^2,\]
where the last inequality follows from \citet[Lemma 5.4]{tu2015low}. This completes the proof of the lemma.

\subsection{Proof of \lemref{RPCA}}\label{app:RPCA_proofs}

Recalling the constrained least squares problem \eqnref{RPCA_program} for the robust PCA problem, we verify that under the conditions of \lemref{RPCA}, the loss function satisfies the assumptions of \thmref{altmin_contract}, i.e. \assumpsref{RSC_joint},~\assumpssref{RSM},~\assumpssref{crossprod}, 
and~\assumpssref{initialization}, with parameters specified below. Before proceeding, observe that for all $Y\in\Yset_0$, we have
\begin{equation}\label{eqn:ell1_bound}
\norm{Y-\Yh}_1  \leq  \norm{Y-\Ys}_1 + \norm{\Yh-\Ys}_1 \leq 2\sqrt{sd}\cdot \fronorm{Y-\Yh} + 4\sqrt{sd}\cdot \fronorm{\Yh-\Ys},
\end{equation}
where the last step holds by the triangle inequality and the fact that $\Ys$ is $sd$-sparse by our assumption.

We now establish the joint restricted strong convexity, restricted smoothness, and initialization conditions, for the least squares loss $\loss(X,Y)=\frac{1}{2}\fronorm{Z-\Aoper(X+Y)}^2$. Note that the cross-product condition (\assumpref{crossprod}) trivially holds with $\mu_x=\mu_y=0$, since the Hessian  $\nabla^2_{XY}\loss(X,Y)$ is constant over all $(X,Y)$. We  use the shorthand $\sigma_r=\sigma_r(\Xs)$ to denote the smallest singular value of $\Xs$.

\paragraph{(Joint RSC.)}

Take $X\in\Xset_0$, $Y\in\Yset_0$. By~\assumpref{RE} (Restricted Eigenvalue), we have that for $\Delta_X=X-\Xh$, $\Delta_Y=Y-\Yh$,
\begin{multline*}
\biginner{\left(\begin{array}{c}\Delta_X\\ \Delta_Y \end{array}\right)}
{\nabla\loss(X,Y)- \nabla\loss(\Xh,\Yh)} = \fronorm{\Aoper(\Delta_X+\Delta_Y)}^2   \\ \geq \alpha_A( \fronorm{\Delta_X}^2 + \fronorm{\Delta_Y}^2 )-  \tau\left(\underbrace{\frac{\log d}{n^2}\norm{\Delta_Y}_1^2 + \sqrt{\frac{d^2\log d}{n^2}}\norm{\Delta_X}_\infty \norm{\Delta_Y}_1}_{(\textnormal{Term 1})} \right).
\end{multline*}

Using the inequality \eqnref{ell1_bound}, and the spikiness constraint, i.e. $\norm{X}_\infty\leq \frac{\spike}{d}$ for $X\in\Xset$, (Term 1) above  is bounded by
\begin{align*}
\textnormal{(Term 1)} &\leq  \frac{4sd\log d}{n^2}\left(  \fronorm{\Delta_Y} + 2\fronorm{\Yh-\Ys} \right)^2 + 4\spike\sqrt{\frac{sd \log d}{n^2}} \left(  \fronorm{\Delta_Y} + 2\fronorm{\Yh-\Ys} \right) 
\\ &\leq  \left[\frac{8sd\log d}{n^2} \fronorm{\Delta_Y}^2 + \frac{32sd\log d}{n^2}\fronorm{\Yh-\Ys}^2\right] +\left[\frac{\alpha_A}{4}  \fronorm{\Delta_Y}^2 + \alpha_A\fronorm{\Yh-\Ys}^2 +  \frac{32\spike^2}{\alpha_A}\frac{sd\log d}{n^2}\right],
\end{align*}
where the second step uses the identity $ab\leq\frac{ca^2}{2} + \frac{b^2}{2c}$ for any $c>0$. Substituting to the inequality above, and using the fact that $\tau\frac{32 sd\log d}{n^2} \leq \alpha_A$, then we have
\[ \biginner{\left(\begin{array}{c}\Delta_X\\ \Delta_Y \end{array}\right)}
{\nabla\loss(X,Y)- \nabla\loss(\Xh,\Yh)} \geq \underbrace{\alpha_A}_{=\alpha_X}\fronorm{\Delta_X}^2 +\underbrace{\frac{\alpha_A}{2}}_{=\alpha_Y}\left[ \fronorm{\Delta_Y}^2 \underbrace{-4\fronorm{\Yh-\Ys}^2 - \frac{64\spike^2}{\alpha_A^2}\frac{sd\log d}{n^2}}_{=\eps_Y^2} \right]. \]

\paragraph{(RSM.)} Again by~\assumpref{RE} (Restricted Eigenvalue), we have that
\[\inner{\Delta_X}{\nabla_X\loss(X,\Yh)-\nabla_X\loss(\Xh,\Yh)} = \fronorm{\Aoper(\Delta_X)}^2 \leq \underbrace{ \beta_A}_{=\beta_X}\fronorm{\Delta_X}^2, \]
and  
\begin{align*}\inner{\Delta_Y}{\nabla_Y\loss(\Xh,Y) - \nabla_Y\loss(\Xh,\Yh)} &= \fronorm{\Aoper(\Delta_Y)}^2 \\&\leq \underbrace{\frac{3\beta_A}{2}}_{=\beta_Y}\fronorm{\Delta_Y}^2 + \frac{\alpha_A}{2}\left(4\fronorm{\Yh-\Ys}^2  + \frac{64\spike^2}{\alpha_A^2}\frac{sd\log d}{n^2} \right).\end{align*}

\paragraph{(Initialization condition.)}

Since $\Yset$ is convex, the initialization condition is trivial for the set $\Yset_0$. For $\Xset_0$, we first bound $\opnorm{\nabla_X\loss(X,Y)}$ for $X\in\Xset_0$ and $Y\in\Yset_0$. Given the  model $Z=\Aoper(\Xs+\Ys)+W$, we have the decomposition 
\[\opnorm{\nabla_X\loss(X,Y)} \leq  \underbrace{\opnorm{\Aoper^*\Aoper(X-\Xs)}}_{\text{(Term 1)}} +\underbrace{\opnorm{\Aoper^*\Aoper(Y-\Ys)}}_{\text{(Term 2)}} + {\opnorm{\Aoper^*(W)}}. \]

Note that $\opnorm{\Aoper^*\Aoper(X-\Xs)}=\inner{\Aoper(X')}{\Aoper(X-\Xs)}$ for some $X'$ with $\rank(X')=1$ and $\fronorm{X'}\leq 1$. By~\assumpref{RE} (Restricted Eigenvalue), then,  $\fronorm{\Aoper(X')}^2 \leq \beta_A$, and we also have $\fronorm{\Aoper(X-\Xs)}^2\leq \beta_A\fronorm{X-\Xs}^2$. Then
\begin{align*}
\text{(Term 1)} = \fronorm{\Aoper(X')}\fronorm{\Aoper(X-\Xs)} \leq \beta_A\fronorm{X-\Xs} \leq 2\beta_A \rho_X,
\end{align*}
where the last inequality holds since $X,\Xs\in\ball_2(X_0,\rho_X)$. Also, by \assumpref{RE}, we have the bound $\fronorm{\Aoper(Y-\Ys)}^2\leq \beta_A\fronorm{Y-\Ys}^2 + \tau\frac{4 sd \log d}{n^2}\fronorm{Y-\Ys}^2\leq \frac{9\beta_A}{8}\fronorm{Y-\Ys}^2$, and so for some $X''$ with $\rank(X'')=1$ and $\fronorm{X''}\leq 1$,
\[\text{(Term 2)} = \inner{\Aoper(X'')}{\Aoper(Y-\Ys)}\leq  \frac{3\sqrt{2}}{4}\beta_A\fronorm{Y-\Ys} \leq 3\beta_A \rho_Y. \]
Putting these bounds together, we have $\opnorm{\nabla_X\loss(X,Y)} \leq 3\beta_A (\rho_X + \rho_Y) + \opnorm{\Aoper^*(W)}$. Now, by \lemref{matrix_lcc}, we know $\gamma_X(\Xset)\leq \frac{5}{4\sigma_r(X)}$, and so
\begin{equation}\label{eqn:matrix_lcc_bound}
\max_{X\in\Xset_0}\gamma_X(\Xset) \leq \frac{5}{4\sigma_r - 8\rho_X} \leq \frac{5}{2\sigma_r} ,\end{equation}
where the first inequality applies Weyl's inequality, while the second inequality uses $\rho_X\leq\frac{1}{4}\sigma_r$. Recalling $\rho_X,\rho_Y\leq c_0\cdot \sigma_r\kappa^{-1}(\Aoper)$ for some sufficiently small $c_0>0$, this implies that the conditions of \lemref{LCC_rho} hold, i.e. $\rho_X < \frac{1}{2\max_{X\in\Xset_0}\gamma_X(\Xset)}$, and in particular, we have $\gamma(\Xset_0)\leq \frac{5}{2\sigma_r}$. Now combining all the pieces, then,
\begin{align*}
2\gamma(\Xset_0)\cdot \left(\opnorm{\nabla_X\loss(\Xh,\Yh)} + \max_{Y\in\Yset_0}\opnorm{\nabla_Y\loss(X_Y,Y)} \right) &\leq 
4 \gamma(\Xset_0)\cdot \max_{X\in\Xset_0,Y\in\Yset_0}\opnorm{\nabla_X\loss(X,Y)} \\ &\leq
\frac{10}{\sigma_r}\cdot\left(3\beta_A\rho_X + 3\beta_A\rho_Y + \opnorm{\Aoper^*(W)} \right) \leq\alpha_A,\end{align*}
where we use $\opnorm{\Aoper^*(W)} \leq \sigma_r\cdot \frac{\alpha_A}{30}$ in the last step. This establishes the initialization condition.

Now by specializing~\thmref{altmin_contract} to the robust PCA problem \eqnref{RPCA_program}, the result of \lemref{RPCA} immediately follows.

\subsection{Proof of \lemref{factor_model}}\label{app:factor_proofs}

Next we turn to prove our claims for the Gaussian factor model, as presented in \eqnref{factor_program}. First, with some algebra, we have the following expression for the gradient and Hessian of $\loss(X,Y)$: for all $\Delta_X,\Delta_{Y}\in\R^{d\times d}$, 
\[\biginner{\left(\begin{array}{c}
	\Delta_X \\ \Delta_Y
	\end{array}\right)}{\nabla\loss(X,Y)} = \textnormal{tr}((\Delta_X+\Delta_Y)^\top (X+Y)^{-1}(X + Y - S_n) (X + Y)^{-1}),\]
and 
\begin{multline*}
\left(\begin{array}{c}
\Delta_X \\ \Delta_Y
\end{array}\right)^\top \nabla^2\loss(X,Y) \left(\begin{array}{c}
\Delta_{X} \\ \Delta_{Y}
\end{array}\right) = \vect{\Delta_X}^\top \mathcal{H}(X,Y) \vect{\Delta_{X}} + \vect{\Delta_Y}^\top \mathcal{H}(X,Y) \vect{\Delta_{Y}} \\
+ 2\vect{\Delta_X}^\top \mathcal{H}(X,Y) \vect{\Delta_{Y}},
\end{multline*}
where $\mathcal{H}(X,Y)$ is a $d^2$-by-$d^2$ matrix, given by
\begin{multline*}\mathcal{H}(X,Y) =  \underbrace{\frac{1}{2}(X+Y)^{-1}(2 S_n - (X+Y)) (X+Y)^{-1} \otimes (X+Y)^{-1}}_{\coloneqq \mathcal{H}_1(X,Y) } \\+ \underbrace{\frac{1}{2} (X+Y)^{-1}\otimes (X+Y)^{-1}(2 S_n - (X+Y)) (X+Y)^{-1}}_{\coloneqq \mathcal{H}_2(X,Y)}.\end{multline*}

In the proof, the following concentration inequality will be used: since $z_i\iidsim \mathcal{N}(0,\Sigmas)$ and $S_n$ is a sample covariance matrix formed by $\{z_i\}_{i=1}^{n}$, with probability at least $1-2e^{-d}$, we have
\begin{align}\label{eqn:factor_concentration}
\opnorm{S_n - \Sigmas} \leq \opnorm{\Sigmas}\opnorm{\Sigmas{}^{-1/2}S_n \Sigmas{}^{-1/2} - \ident_d} 
\leq 3\lammax(\Sigmas)\sqrt{\frac{d}{n}},
\end{align}
where the second step holds by a concentration bound on the extreme singular values of a standard Gaussian ensemble \citet{davidson2001local}. 

We calculate a few inequalities to use later. Recall $\rho_X,\rho_Y\leq c_0\cdot \min\{ \sigma_r(\Xs)\kappa^{-3}(\Sigmas),\lammin(\Sigmas)\kappa^{-4}(\Sigmas)\}$ for a sufficiently small $c_0>0$.  For $X\in\ball_2(X_0,\rho_X)$ and $Y\in\ball_2(Y_0,\rho_Y)$, since $\Xs\in\Xset_0, \Ys\in\Yset_0$ while $\Sigmas=\Xs+\Ys$, 
\[\opnorm{X+Y-\Sigmas} \leq \opnorm{X+Y-X_0-Y_0} +\opnorm{X_0+Y_0-\Sigmas} \leq 2\rho_X + 2\rho_Y \leq \frac{\lammin(\Sigmas)}{4},\]
where the last inequality holds since $\rho_X,\rho_Y\leq \frac{\lammin(\Sigmas)}{16}$.
Applying Weyl's inequality, this yields
\begin{align}\label{eqn:factor_eig_bound1}
\frac{3}{4}\lammin(\Sigmas) \leq \lammin(X + Y) \leq \lammax(X + Y) \leq \frac{5}{4}\lammax(\Sigmas).
\end{align}
Applying Weyl's inequality again, and using the inequality \eqnref{factor_concentration}, we also have
\begin{align}\label{eqn:factor_eig_bound2}
\frac{1}{2}\lammin(\Sigmas) \leq \lammin(2S_n - X-Y) \leq \lammax(2S_n -X-Y) \leq   \frac{3}{2}\lammax(\Sigmas),
\end{align}
where we use the assumption $\sqrt{\frac{d}{n}}\leq \frac{\kappa^{-1}(\Sigmas)}{24}$. In particular, combining these bounds, and using  standard properties of the Kronecker product, we have that
\begin{align}\label{eqn:factor_eig_bound3}
\frac{32}{125}\frac{\kappa^{-1}(\Sigmas)}{\lammax^2(\Sigmas)}\leq \lammin(\mathcal{H}(X,Y)) \leq \lammax(\mathcal{H}(X,Y)) \leq \frac{32}{9}\frac{\kappa(\Sigmas)}{\lammin^2(\Sigmas)}.
\end{align}

Finally, due to the spikiness, $\norm{X}_\infty\leq \frac{\spike}{d}$, and the $\ell_1$ norm inequality \eqnref{ell1_bound}, we have the following finite bound on the inner product between the low-rank and sparse components: for all $X\in\Xset_0$, all $Y\in\Yset_0$, writing $\Delta_X=X-\Xh, \Delta_Y = Y-\Yh$,
\begin{equation}\label{eqn:inner_prod}
\inner{\Delta_X}{\Delta_Y} \leq \norm{\Delta_X}_\infty \norm{\Delta_Y}_1  \leq  4\spike\sqrt{\frac{s}{d}}\cdot \fronorm{\Delta_Y} + 8\spike\sqrt{\frac{s}{d}}\cdot \fronorm{\Yh-\Ys}.
\end{equation}

 Throughout the proof, we use the shorthand notation $\sigma_r=\sigma_r(\Xs)$. 

\paragraph{(Joint RSC.)} Write $\Delta_X=X-\Xh, \Delta_Y = Y-\Yh$, then by Taylor's theorem, it is sufficient to lower bound
\begin{align*}
\left(\begin{array}{c}
\Delta_X \\ \Delta_Y
\end{array}\right)^\top \nabla^2\loss(X(t),Y(t)) \left(\begin{array}{c}
{\Delta_X} \\ \Delta_Y
\end{array}\right)  
=\left(\vect{\Delta_X} + \vect{\Delta_Y}\right)^\top  \mathcal{H}(X(t),Y(t))\left(\vect{\Delta_X} + \vect{\Delta_Y}\right),
\end{align*}
where, for some $t\in [0,1]$, $X(t)=(1-t)X + t\Xh$ and $Y(t)=(1-t)Y+t\Yh$. By \eqnref{factor_eig_bound3}, the right-hand side is lower bounded by $ \lammin(\mathcal{H}(X(t),Y(t)))\cdot \norm{\vect{\Delta_X} + \vect{\Delta_Y}}_2^2$. We also have
\begin{align*}\Norm{\vect{\Delta_X} + \vect{\Delta_Y}}_2^2 
&\geq \fronorm{\Delta_X}^2 + \fronorm{\Delta_Y}^2 -  8\spike\sqrt{\frac{s}{d}}\cdot \fronorm{\Delta_Y} - 16\spike\sqrt{\frac{s}{d}}\cdot \fronorm{\Yh-\Ys} \\
&\geq \fronorm{\Delta_X}^2 +\frac{1}{2} \fronorm{\Delta_Y}^2 - 16\fronorm{\Yh-\Ys}^2 - \spike^2\frac{36s}{d}, \end{align*}
where the first step applies~\eqnref{inner_prod}, while the second step uses the inequality $ab\leq \frac{ca^2}{2} + \frac{b^2}{2c}$. We thus have
\begin{equation*}
\left(\begin{array}{c}
\Delta_X \\ \Delta_Y
\end{array}\right)^\top \nabla^2\loss(X(t),Y(t)) \left(\begin{array}{c}
{\Delta_X} \\ \Delta_Y
\end{array}\right)  \geq  \underbrace{  \frac{32}{125}\frac{\kappa^{-1}(\Sigmas)}{\lammax^2(\Sigmas)}}_{=\alpha_X,2\cdot \alpha_Y} \left(\fronorm{\Delta_X}^2 +\frac{1}{2}\fronorm{\Delta_Y}^2 - \underbrace{16\fronorm{\Yh-\Ys}^2  - \spike^2 \frac{36s}{d}}_{=\eps_Y^2} \right).
\end{equation*}

\paragraph{(RSM.)} By  Taylor's theorem, and using the inequality \eqnref{factor_eig_bound3}, it is easy to see that
\[\inner{\Delta_X}{\nabla_X\loss(X,\Yh) - \nabla_X\loss(\Xh,\Yh) }  \leq \underbrace{\frac{32}{9}\frac{\kappa(\Sigmas)}{\lammin^2(\Sigmas)}}_{=\beta_Y} \fronorm{\Delta_X}^2,\]
and analogously for the $Y$ variable.

\paragraph{(Cross-product bound.)} 
As discussed in \secref{assumptions} following the \assumpref{crossprod}, in order to establish the cross-product condition, it suffices to bound $\opnorm{\nabla^2_{XY}\loss(X,Y(t)) - \nabla^2_{XY}\loss(X(t'),Y)} $, where $X(t')=(1-t')X+t'\Xh$ and $Y(t)=(1-t)Y+t\Yh$. Here we  only focus on bounding the term $ \opnorm{\nabla^2_{XY}\loss(X,Y(t)) - \nabla^2_{XY}\loss(X,Y)}$; by symmetry, a similar bound holds for  $ \opnorm{\nabla^2_{XY}\loss(X(t'),Y) - \nabla^2_{XY}\loss(X,Y)}$, which, combined with the triangle inequality, gives the desired bound.

Furthermore, by the property of Kronecker product, the operator norms of $\mathcal{H}_1$ and $\mathcal{H}_2$ are equal, and so
\begin{align*}
\opnorm{\nabla^2_{XY}\loss(X,Y(t)) - \nabla^2_{XY}\loss(X,Y)} &= \opnorm{\mathcal{H}(X,Y(t))-\mathcal{H}(X,Y)} \leq 2\opnorm{\mathcal{H}_1(X,Y(t))-\mathcal{H}_1(X,Y)},
\end{align*}
where the inequality holds since $\mathcal{H}=\mathcal{H}_1+\mathcal{H}_2$ and applying the triangle inequality; therefore, it suffices to have a bound on the term $\opnorm{\mathcal{H}_1(X,Y(t))-\mathcal{H}_1(X,Y)}$.

Now let $\Delta\mathcal{H}_{1}= (X+Y(t))^{-1}(2 S_n - (X+Y(t))) (X+Y(t))^{-1} - (X+Y)^{-1}(2 S_n - (X+Y)) (X+Y)^{-1}$. Simple calculation shows that
\begin{multline*}
\mathcal{H}_1(X,Y(t))-\mathcal{H}_1(X,Y) = \frac{1}{2}\Delta\mathcal{H}_{1}\otimes (X+Y(t))^{-1} \\ + \frac{1}{2} (X+Y)^{-1}(2 S_n - (X+Y)) (X+Y)^{-1} \otimes \left( (X+Y(t))^{-1} - (X+Y)^{-1} \right).
\end{multline*}
$\Delta\mathcal{H}_{1}$ is further decomposed into the sum
\begin{multline*}
\Delta\mathcal{H}_{1} = \left((X+Y(t))^{-1} - (X+Y)^{-1} \right)(2 S_n - (X+Y(t))) (X+Y(t))^{-1} \\+ (X+Y)^{-1} (Y - Y(t)) (X+Y(t))^{-1} + (X+Y)^{-1}(2 S_n - (X+Y(t))) \left((X+Y(t))^{-1} - (X+Y)^{-1} \right).
\end{multline*}
Meanwhile, by the inequalities \eqnref{factor_eig_bound1} and \eqnref{factor_eig_bound2}, we have that
\[\opnorm{(X+Y(t))^{-1}}, \opnorm{(X+Y)^{-1}} \leq \frac{4}{3\lammin(\Sigmas)}\text{ and }\opnorm{2S_n - (X+Y(t))}, \opnorm{2S_n - (X+Y)}\leq \frac{3\lammax(\Sigmas)}{2}.\]
Using the identity $A^{-1}+B^{-1}=A^{-1}(A+B)B^{-1}$, we can also see that
\[\opnorm{(X+Y(t))^{-1}- (X+Y)^{-1}}\leq \frac{16}{9\lammin^2(\Sigmas)}\cdot t\opnorm{\Delta_Y},\]
so combining these inequalities, we obtain
\[\opnorm{\Delta\mathcal{H}_1} \leq \left( \frac{64}{9}\frac{\lammax(\Sigmas)}{\lammin^3(\Sigmas)}\right)\cdot \opnorm{\Delta_Y} +\left(\frac{16}{9}\frac{1}{\lammin^2(\Sigmas)}\right)\cdot \opnorm{\Delta_Y}.\]
Then:
\begin{align*}
&\opnorm{\mathcal{H}_1(X,Y(t))-\mathcal{H}_1(X,Y)} \\& \leq \frac{1}{2}\opnorm{\Delta\mathcal{H}_1}\opnorm{(X+Y(t))^{-1}} + \frac{1}{2}\opnorm{ (X+Y)^{-1}(2 S_n - (X+Y)) (X+Y)^{-1} } \opnorm{(X+Y(t))^{-1}- (X+Y)^{-1}}\\
&\leq \left( \frac{192}{27}\frac{\lammax(\Sigmas)}{\lammin^4(\Sigmas)}\right) \cdot \opnorm{\Delta_Y} + \left(\frac{32}{27}\frac{1}{\lammin^3(\Sigmas)}\right)\cdot \opnorm{\Delta_Y}  \leq \left( \frac{224}{27}\frac{\lammax(\Sigmas)}{\lammin^4(\Sigmas)}\right)\cdot \opnorm{\Delta_Y}.
\end{align*}
Returning to the cross product condition, this implies that,
\[\opnorm{\nabla^2_{XY}\loss(X,Y(t)) - \nabla^2_{XY}\loss(X,Y)}\leq \left(\frac{448}{27}\frac{\lammax(\Sigmas)}{\lammin^4(\Sigmas)}\right)\cdot \opnorm{\Delta_Y}, \]
and in particular, by symmetry, we have
\[\opnorm{\nabla^2_{XY}\loss(X,Y(t)) - \nabla^2_{XY}\loss(X(t'),Y)} \leq \left(\frac{448}{27}\frac{\lammax(\Sigmas)}{\lammin^4(\Sigmas)} \right)\cdot \left(\opnorm{\Delta_X} + \opnorm{\Delta_Y}\right).\]
To summarize, we have shown that $\mu_X=\mu_Y = \left[\frac{896}{27}\frac{\lammax(\Sigmas)}{\lammin^4(\Sigmas)}\right] \cdot \left(\rho_X + \rho_Y\right)$. By choosing $c_0$ sufficiently small, this gives the claim $\mu_X=\mu_Y\leq \frac{16}{125}\frac{\lammin(\Sigmas)}{\lammax^3(\Sigmas)}$ as desired.

\paragraph{(Initialization condition.)}  

To prove the initialization condition, it is sufficient to bound the quantity $4\gamma(\Xset_0)\cdot\max_{X\in\Xset_0,Y\in\Yset_0}\opnorm{\nabla_X\loss(X,Y)}$. Note that for any $X\in\Xset_0,Y\in\Yset_0$,
\begin{align*}
\opnorm{\nabla_X\loss(X,Y)} &= \opnorm{(X+Y)^{-1}(X+Y-S_n)(X+Y)^{-1}} \\ &\leq \frac{16}{9\lammin^2(\Sigmas)}\cdot \left(\opnorm{X+Y-\Sigmas} + \opnorm{S_n - \Sigmas}\right)
\leq \left( \frac{8}{625}\sigma_r\right)\cdot \frac{\lammin(\Sigmas)}{\lammax^3(\Sigmas)},
\end{align*}
where the first inequality applies~\eqnref{factor_eig_bound1}, while the second inequality applies the bound $\opnorm{X+Y-\Sigmas}\leq 2(\rho_X+\rho_Y)$, and the concentration bound~\eqnref{factor_concentration}, as well as our assumptions on the radii and the sample size \eqnref{factor_sample} (where we choose $c_0,c_1>0$ to be sufficiently small). Also, by the same reasoning to the equation \eqnref{matrix_lcc_bound}, we  have $\gamma(\Xset_0)\leq \frac{5}{2\sigma_r}$. Therefore,
\[4\gamma(\Xset_0)\cdot\max_{X\in\Xset_0,Y\in\Yset_0}\opnorm{\nabla_X\loss(X,Y)} \leq \frac{10}{\sigma_r}\cdot \left(\frac{8}{625}\sigma_r\right)\frac{\lammin(\Sigmas)}{\lammax^3(\Sigmas)} = \frac{16}{125}\frac{\lammin(\Sigmas)}{\lammax^3(\Sigmas)} \leq \alpha_X-\mu_X,\]
as desired.

 \subsection{Proof of \lemref{multitask_regression}}
 
 Recall the expression for the negative log-likelihood function 
 \[\loss(X,\Theta) = -\log\det(\Theta) + \frac{1}{n}\sum_{i=1}^{n} (z_i - X\phi_i)^\top \Theta (z_i- X\phi_i).\]

Then we can calculate
\[\nabla_X\loss(X,\Theta) = \frac{2}{n}\sum_{i=1}^{n} \Theta(X\phi_i - z_i) \phi_i^\top \; \text{ and } \; \nabla_\Theta\loss(X,\Theta) =  -\Theta^{-1} + \frac{1}{n}\sum_{i=1}^{n}(z_i - X\phi_i) (z_i - X\phi_i)^\top,\]
and
\begin{align*}
\inner{\Delta_X}{\nabla^2_{XX}\loss(X,\Theta)\Delta_{X}} &= \frac{2}{n}\sum_{i=1}^{n} \phi_i^\top \Delta_X^\top \Theta  \Delta_{X} \phi_i,\\
\inner{\Delta_X}{\nabla^2_{X\Theta}\loss(X,\Theta)\Delta_{\Theta}} &= \frac{2}{n}\sum_{i=1}^{n} \phi_i^\top \Delta_X^\top \Delta_\Theta(X \phi_i -z_i),\\
 \inner{\Delta_{\Theta}}{\nabla^2_{\Theta\Theta}\loss(X,\Theta)\Delta_{\Theta}} &= \vect{\Delta_{\Theta}}^\top\left(\Theta^{-1}\otimes \Theta^{-1}\right)\vect{\Delta_{\Theta}}.
\end{align*}

Throughout we use the shorthand notation $\sigma_r=\sigma_r(\Xs)$. Recall that the radii are chosen 
to satisfy $\rho_X\leq c_0\cdot \sigma_r\kappa^{-1}(\Thetas)\kappa^{-1}(\Sigma_\phi)$ and $\rho_\Theta\leq c_0\cdot \lammin(\Thetas)\kappa^{-1}(\Sigma_\phi)$ for some small $c_0>0$. Then, according to Weyl's inequality, for any $\Theta\in\Thetaset_0$, its minimum and maximum eigenvalues are bounded by
\begin{equation}\label{eqn:multitask_eig_bound}
\frac{\lammin(\Thetas)}{2} \leq \lammin(\Theta)\leq \lammax(\Theta) \leq \frac{3\lammax(\Thetas)}{2},
\end{equation}
where we use $\fronorm{\Theta-\Thetas}\leq 2\rho_\Theta$, and $\rho_\Theta \leq \frac{\lammin(\Thetas)}{4}$.

We will use the following two concentration results: first, by~\citet[Lemma 2]{negahban2011estimation}, with probability at least $1-4\exp(-n/2)$, we have the following bounds of the form:
\begin{equation}\label{eqn:multitask_concen1} \lammin\left(\frac{1}{n}\sum_{i=1}^{n}\phi_i\phi_i^\top \right)\geq \frac{\lammin(\Sigma_\phi)}{9} \;\text{ and }\;  \lammax\left(\frac{1}{n}\sum_{i=1}^{n}\phi_i\phi_i^\top\right) \leq 9\lammax(\Sigma_\phi).\end{equation}

Next, letting $\widetilde{\epsilon}_i\iidsim\textnormal{N}(0,\ident_m)$, it has been shown in \citet[Lemma 3]{negahban2011estimation} that for some $c,c'>0$, with probability at least $1-c\exp(-c'(m+d))$,
\begin{equation}\label{eqn:multitask_concen2}\opNorm{\frac{1}{n}\sum_{i=1}^{n}\widetilde{\epsilon}_i \phi_i^\top} \leq 5\sqrt{\lammax(\Sigma_\phi)}\sqrt{\frac{m+d}{n}}. \end{equation}

Now, we turn to verifying \lemref{multitask_regression}:

\paragraph{(Joint RSC.)}

Take $X\in\Xset_0$, $\Theta\in\Thetaset_0$. By Taylor's theorem, we have  ($\Delta_X=X-\Xh$, $\Delta_\Theta=\Theta-\Thetah$)
\[\biginner{\left(\begin{array}{c}
\Delta_X\\ \Delta_\Theta
\end{array}\right)}{\nabla\loss(X,\Theta)-\nabla\loss(\Xh,\Thetah)} = \left(\begin{array}{c}
\Delta_X \\ 
\Delta_\Theta
\end{array} \right)^\top \nabla^2 \loss(X(t),\Theta(t)) \left(\begin{array}{c}
\Delta_X \\ 
\Delta_\Theta
\end{array} \right), \]
where we write $X(t) = (1-t)X + t\Xh$ and $\Theta(t)=(1-t)\Theta + t\Thetah$ for some $t\in(0,1)$. Using the expression for the Hessian operator, and substituting the observational model $z_i = \Xs \phi_i + \epsilon$, we have the following decomposition:
\begin{multline*}
 \left(\begin{array}{c}
\Delta_X \\ 
\Delta_\Theta
\end{array} \right)^\top \nabla^2 \loss(X(t),\Theta(t)) \left(\begin{array}{c}
\Delta_X \\ 
\Delta_\Theta
\end{array} \right)= \underbrace{ \frac{2}{n}\sum_{i=1}^{n} \phi_i^\top \Delta_X^\top \Theta(t) \Delta_X\phi_i}_{\textnormal{(Term 1)}} \\
+ \underbrace{ \frac{2}{n}\sum_{i=1}^{n}   \phi_i^\top \Delta_X^\top \Delta_\Theta  (X(t)-\Xs) \phi_i}_{\textnormal{(Term 2)}} - \underbrace{ \frac{2}{n}\sum_{i=1}^{n}  \phi_i^\top \Delta_X^\top \Delta_\Theta\cdot \epsilon_i}_{\textnormal{(Term 3)}} \\ +\underbrace{  \vect{\Delta_\Theta}^\top \left(\Theta(t)^{-1}\otimes \Theta(t)^{-1} \right) \vect{ \Delta_\Theta}}_{\textnormal{(Term 4)}}.
\end{multline*}

(Term 1) is lower bounded by
\[\text{(Term 1)} \geq 2\lammin\left(\frac{1}{n}\sum_{i=1}^{n}\phi_i \phi_i^\top\right)\cdot \lammin(\Theta(t)) \fronorm{\Delta_X}^2\geq \frac{\lammin(\Thetas)\lammin(\Sigma_\phi)}{9}\fronorm{\Delta_X}^2,\]
where the second step uses the inequalities \eqnref{multitask_eig_bound} and \eqnref{multitask_concen1}. For (Term 2), we further decompose it as
\[\text{(Term 2)} = (1-t)\cdot \frac{2}{n}\sum_{i=1}^{n}  \phi_i^\top \Delta_X^\top \Delta_\Theta  \Delta_X\phi_i +t\cdot \frac{2}{n}\sum_{i=1}^{n}   \phi_i^\top \Delta_X^\top \Delta_\Theta   (\Xh-\Xs)\phi_i.\]
Then the first term is bounded by 
\[4\rho_\Theta \cdot \lammax\left(\frac{1}{n}\sum_{i=1}^{n} \phi_i \phi_i^\top\right)\fronorm{\Delta_X}^2 \leq  \frac{\lammin(\Thetas) \lammin(\Sigma_\phi)}{54}\fronorm{\Delta_X}^2,\]
where the inequality uses the bound on the radius $\rho_{\Theta}$ (by choosing $c_0\leq \frac{1}{36\cdot 54}$), and the concentration bound~\eqnref{multitask_concen1}, while the second part of (Term 2) is bounded by
\begin{align*}
 \frac{2}{n}\sum_{i=1}^{n}   \phi_i^\top \Delta_X^\top \Delta_\Theta   (\Xh-\Xs)\phi_i &\leq \frac{4\rho_\Theta}{n}\sum_{i=1}^{n}  \norm{\Delta_X\phi_i}_2 \norm{(\Xh-\Xs)\phi_i}_2 \\ &\leq   \frac{2\rho_\Theta}{n}\sum_{i=1}^{n}  \norm{\Delta_X \phi_i}_2^2 +  \frac{2\rho_\Theta}{n}\sum_{i=1}^{n} \norm{(\Xh-\Xs)\phi_i}_2^2 \\ &\leq \frac{\lammin(\Thetas) \lammin(\Sigma_\phi)}{108}\fronorm{\Delta_X}^2 + \frac{\lammin(\Thetas) \lammin(\Sigma_\phi)}{108}\fronorm{\Xh-\Xs}^2,
\end{align*}
where the first step applies the Cauchy-Schwarz inequality; the second step applies the identity $2ab\leq a^2 + b^2$; and the third uses $\rho_\Theta\leq c_0\cdot\lammin(\Thetas)\kappa^{-1}(\Sigma_\phi)$, and applies~\eqnref{multitask_concen1}.
Combining the two, then,
\[\text{(Term 2)} \leq  \frac{\lammin(\Thetas) \lammin(\Sigma_\phi)}{36}\fronorm{\Delta_X}^2 + \frac{\lammin(\Thetas) \lammin(\Sigma_\phi)}{108}\fronorm{\Xh-\Xs}^2.\]
Next, using the inequality $\inner{a}{b}\leq\nucnorm{a}\opnorm{b}$, we find that 
\begin{multline*}
\text{(Term 3)} \leq 2\nucnorm{\Delta_X}\opNorm{\frac{1}{n}\sum_{i=1}^{n}\Delta_\Theta\cdot \epsilon_i \phi_i^\top} \leq \frac{\rho_\Theta}{\sqrt{\lammin(\Thetas)}}\cdot 2\sqrt{2r}\fronorm{\Delta_X}\opNorm{\frac{1}{n}\sum_{i=1}^{n}\widetilde{\epsilon}_i \phi_i^\top} \leq \\
 \frac{\rho_\Theta\sqrt{\lammax(\Sigma_\phi)}}{\sqrt{\lammin(\Thetas)}}\cdot 10\sqrt{2r}\fronorm{\Delta_X} \sqrt{\frac{m+d}{n}},
\end{multline*}
where the second step follows since $X-\Xh$ is of rank $2r$, and $\epsilon_i=(\Thetas)^{-1/2}\cdot \widetilde{\epsilon}_i$ for $\widetilde{\epsilon}_i\iidsim \mathcal{N}(0,\ident_m)$, and the next step applies the concentration bound \eqnref{multitask_concen2}. Using the identity $ab\leq \frac{ca^2}{2} + \frac{b^2}{2c}$, and the bound on $\rho_\Theta$, then
\[\text{(Term 3)} \leq \frac{\lammin(\Thetas) \lammin(\Sigma_\phi)}{36}\fronorm{\Delta_X}^2 + \left( \frac{25}{13122}\cdot \frac{r(m+d)}{n}\right) \frac{\lammin(\Sigma_\phi)}{\lammax(\Sigma_\phi)}. \]

 Lastly, by \eqnref{multitask_eig_bound}, the minimum eigenvalue of $\Theta(t)^{-1}$ is lower bounded by $\frac{2}{3\lammax(\Thetas)}$, so it follows that
\[\textnormal{(Term 4)} \geq \frac{4}{9\lammax^2(\Thetas)} \fronorm{\Delta_\Theta}^2 .\]
Putting all the bounds together, we have
\begin{multline*}\inner{\left(\begin{array}{c}
\Delta_X\\ \Delta_\Theta
\end{array}\right)}{\nabla\loss(X,\Theta)-\nabla\loss(\Xh,\Thetah)}  \geq \underbrace{\frac{\lammin(\Thetas) \lammin(\Sigma_\phi)}{18} }_{=\alpha_X}\left(\fronorm{\Delta_X}^2 - \frac{1}{6}\fronorm{\Xh-\Xs}^2  \right.\\\left.  - \frac{25}{729} \frac{r(m+d)}{n} \frac{1}{\lammin(\Thetas) \lammax(\Sigma_\phi)}\right)
 +  \underbrace{\frac{4}{9\lammax^2(\Thetas)} }_{=\alpha_\Theta}\fronorm{\Delta_\Theta}^2.\end{multline*}

\paragraph{(RSM.)} 
We can easily see that for $X$,
\begin{align*}
\inner{\Delta_X}{\nabla_X\loss(X,\Thetah) - \nabla_X\loss(\Xh,\Thetah)} &=  \frac{2}{n} \sum_{i=1}^{n}\phi_i^\top \Delta_X^\top \Thetah \Delta_X\phi_i \\ &\leq \underbrace{27\lammax(\Thetas)\lammax(\Sigma_\phi)}_{=\beta_X}\fronorm{\Delta_X}^2,
 \end{align*}
 where the inequality applies~\eqnref{multitask_eig_bound} and \eqnref{multitask_concen1}. 
 Meanwhile, by Taylor's theorem, for some $t\in [0,1]$,
\begin{align*}\inner{\Delta_\Theta}{\nabla_\Theta\loss(\Xh,\Theta)-\nabla_\Theta\loss(\Xh,\Thetah)} &=  \vect{\Delta_\Theta}^\top\nabla^2_{\Theta\Theta}\loss(\Xh,(1-t)\Theta+t\Thetah)\vect{\Delta_\Theta}  \\&\leq  \underbrace{\frac{4}{\lammin^2(\Thetas)} }_{=\beta_\Theta}\fronorm{ \Delta_\Theta}^2,\end{align*}
where the inequality applies~\eqnref{multitask_eig_bound}. This proves the desired results.

\paragraph{(Cross-product bound.)} 
Let $X\in\Xset_0$, $\Theta\in\Thetaset_0$. Then, by Taylor's theorem, for some $t,t'\in[0,1]$,
\begin{align*}
&|\inner{\Delta_X}{\nabla_X\loss(X,\Theta) - \nabla_X\loss(X,\Thetah)} - \inner{\Delta_\Theta}{\nabla_\Theta\loss(X,\Theta)-\nabla_\Theta\loss(\Xh,\Theta)}| \\ &\leq \vect{\Delta_X}^\top \left(\nabla^2_{X\Theta}\loss(X,t\Theta+(1-t)\Thetah) - \nabla^2_{X\Theta}\loss(t'X+(1-t')\Xh,\Theta) \right) \vect{\Delta_\Theta}. \\&= \frac{2(1-t')}{n}\sum_{i=1}^{n}\phi_i^\top \Delta_X^\top \Delta_\Theta \Delta_X \phi_i \leq \underbrace{\frac{\lammin(\Thetas) \lammin(\Sigma_\phi)}{54}}_{=\frac{\mu_X}{2}}\fronorm{\Delta_X}^2,
\end{align*}
where the second step follows from the expression of the Hessian operator $\nabla^2_{X\Theta}\loss$, and the last step applies~\eqnref{multitask_eig_bound} and \eqnref{multitask_concen1}. This proves the cross-product condition, with $\mu_X=\frac{\lammin(\Thetas) \lammin(\Sigma_\phi)}{27}$ and $\mu_\Theta = 0$.

\paragraph{(Initialization condition.)} 

It is shown in~\citep[Lemma 7]{barber2017gradient} that $\gamma_X(\Xset)=\frac{1}{2\sigma_r(X)}$, so we have
\[ \max_{X\in\Xset_0} \gamma_X(\Xset) \leq \frac{1}{2\sigma_r -4\rho_X} \leq \frac{1}{\sigma_r},\]
where the first inequality applies Weyl's inequality, and the next inequality uses $\rho_X\leq \frac{1}{4}\sigma_r$. In particular, this shows that the conditions of \lemref{LCC_rho} is satisfied, and so we have $\gamma(\Xset_0)\leq \frac{1}{\sigma_r}$.

 Next, we bound the gradient term $\opnorm{\nabla_X\loss(X,\Theta)}$. Given the observational model $z_i=\Xs\phi_i +\epsilon_i$, we can decompose the gradient into the two terms
\[\opnorm{\nabla_X\loss(X,\Theta)} \leq \opNorm{\frac{2}{n}\sum_{i=1}^{n}\Theta(X-\Xs)\phi_i\phi_i^\top} + \opNorm{\frac{2}{n}\sum_{i=1}^{n}\Theta \cdot \epsilon_i \phi_i^\top}. \]
Using the inequalities \eqnref{multitask_eig_bound} and \eqnref{multitask_concen1}, the first term is  bounded by $54\rho_X \cdot \lammax(\Thetas) \lammax(\Sigma_\phi)$, whereas we can bound the second term as
\begin{align*}
\opNorm{\frac{2}{n}\sum_{i=1}^{n}\Theta\cdot \epsilon_i \phi_i^\top} \leq \frac{3\lammax(\Thetas)}{\sqrt{\lammin(\Thetas)}}\cdot \opNorm{\frac{1}{n}\sum_{i=1}^{n}\widetilde{\epsilon}_i \phi_i^\top} \leq   \frac{15\lammax(\Thetas)\sqrt{\lammax(\Sigma_\phi)}}{\sqrt{\lammin(\Thetas)}}\sqrt{\frac{m+d}{n}},
\end{align*}
where the steps use the inequalities \eqnref{multitask_eig_bound} and \eqnref{multitask_concen2}.
Combining the two, and using the bound on $\rho_X$, in addition to the assumption \eqnref{multitask_sample}, for sufficiently small $c_0, c_1>0$, 
\[\max_{X\in\Xset_0,\Theta\in\Thetaset_0}\opnorm{\nabla_X\loss(X,\Theta)} \leq \sigma_r\cdot \frac{\lammin(\Thetas)\lammin(\Sigma_\phi)}{216}, \]
and therefore
\[4\gamma(\Xset_0)\cdot\max_{X\in\Xset_0,\Theta\in\Thetaset_0}\opnorm{\nabla_X\loss(X,\Theta)} \leq \frac{\lammin(\Thetas) \lammin(\Sigma_\phi)}{54} = \alpha_X-\mu_X,\]
which completes the proof of \lemref{multitask_regression}.

\end{document}